\documentclass{amsart}
\usepackage{amsmath}
\usepackage{amsfonts}
\usepackage{amssymb}
\usepackage{graphicx}
\newtheorem{theorem}{Theorem}[section]
\newtheorem{lemma}[theorem]{Lemma}
\newtheorem{corollary}[theorem]{Corollary}
\newtheorem{proposition}[theorem]{Proposition}
\newtheorem{conjecture}[theorem]{Conjecture}
\newtheorem{remark}[theorem]{Remark}

\newcommand{\reals}{\mathbb{R}}
\newcommand{\complexes}{\mathbb{C}}
\newcommand{\rationals}{\mathbb{Q}}
\newcommand{\integers}{\mathbb{Z}}

\newcommand{\conv}{\mathrm{conv}\,\,}

\newcommand{\str}{\mathrm{star}}

\newcommand{\assoc}{\mathcal A}

\begin{document}
\title{The signature of a toric variety}
\author{Naichung Conan Leung}
\address{School of Mathematics\\
University of Minnesota\\
Minneapolis, MN 55455 USA}
\email[N.C. Leung]{leung@math.umn.edu}
\author{Victor Reiner}
\address{School of Mathematics\\
University of Minnesota\\
Minneapolis, MN 55455 USA}
\email[V. Reiner]{reiner@math.umn.edu}
\thanks{First, second authors partially supported by NSF grants DMS-9803616 and
DMS-9877047 respectively.}

\begin{abstract}
We identify a combinatorial quantity (the alternating sum of the $h$-vector)
defined for any simple polytope as the signature of a toric variety. This
quantity was introduced by Charney and Davis in their work, which in
particular showed that its non-negativity is closely related to a conjecture
of Hopf on the Euler characteristic of a non-positively curved manifold.

We prove positive (or non-negative) lower bounds for this quantity under
geometric hypotheses on the polytope, and in particular, resolve a special
case of their conjecture. These hypotheses lead to ampleness (or weaker
conditions) for certain line bundles on toric divisors, and then the lower
bounds follow from calculations using the Hirzebruch Signature Formula.

Moreoever, we show that under these hypotheses on the polytope, the $i^{th}$
$L$-class of the corresponding toric variety is $(-1)^{i}$ times an effective
class for any $i$.
\end{abstract}\maketitle

\section{Introduction}

Much attention in combinatorial geometry has centered on the problem of
characterizing which non-negative integer sequences $(f_{0},f_{1},\ldots
,f_{d})$ can be the \textit{$f$-vector} $f(P)$ of a $d$-dimensional convex
polytope $P $, that is, $f_{i}$ is the number of $i$-dimensional faces of $P$;
see \cite{BilleraBjorner} for a nice survey.

For the class of \textit{simple polytopes}, this problem was completely solved
by the combined work of Billera and Lee \cite{BilleraLee} and of Stanley
\cite{Stanley-g-theorem}. A simple $d$-dimensional polytope is one in which
every vertex lies on exactly $d$ edges. McMullen's \textit{$g$-conjecture}
(now the \textit{$g$-theorem}) gives necessary \cite{Stanley-g-theorem} and
sufficient \cite{BilleraLee} conditions for $(f_{0},f_{1},\ldots,f_{d})$ to be
the $f$-vector of a simple $d$-dimensional polytope. Stanley's proof of the
necessity of these conditions showed that they have a very natural phrasing in
terms of the cohomology of the toric variety $X_{\Delta}$ associated to the
\textit{(inner) normal fan} $\Delta$ of $P$, and then the Hard Lefschetz
Theorem for $X_{\Delta}$ played a crucial role. This construction of
$X_{\Delta}$ from $\Delta$ requires that $P$ be \textit{rational}, i.e. that
its vertices all have rational coordinates with respect to some lattice, which
can be achieved by a small perturbation that does not affect $f(P)$. Later,
McMullen \cite{McMullen} demonstrated that one can construct a ring $\Pi(P)$,
isomorphic (with a doubling of the grading) to the cohomology ring of
$X_{\Delta}$ if $P$ is rational, and proved that $\Pi(P)$ formally satisfies
the Hard Lefschetz Theorem, using only tools from convex geometry. In
particular, he recovered the necessity of the conditions of the $g$-theorem in
this way.

This paper shares a similar spirit with Stanley's proof. We attempt to use
further facts about the geometry of $X_{\Delta}$ to deduce information about
the $f$-vector $f(P)$ under certain hypotheses on $P$. The starting point of
our investigation is an interpretation of the alternating sum of the
\textit{$h$-vector} which follows from the Hard Lefschetz Theorem. Recall that
for a simple polytope $P$, the $h$-vector is the sequence $h(P)=(h_{0}%
,h_{1},\ldots,h_{d})$ defined as follows. If we let $f(P,t) := \sum_{i=0}^{d}
f_{i}(P) t^{i} $, then
\[
h(P,t) := \sum_{i=0}^{d} h_{i}(P) t^{i} = f(P,t-1).
\]
The $h$-vector has a topological interpretation: $h_{i}$ is the $2i^{th}$
Betti number for $X_{\Delta}$, or the dimension of the $i^{th}$-graded
component in McMullen's ring $\Pi(P)$. Part of the conditions of the
$g$-theorem are the Dehn-Sommerville equations $h_{i} = h_{d-i}$, which
reflect Poincar\'e duality for $X_{\Delta}$.

Define the alternating sum
\[
\begin{aligned} \sigma(P): & = \sum_{i=0}%
^d (-1)^i h_i(P) \\ [& = h(P,-1) = f(P,-2) = \sum_{i=0}%
^d f_i(P) (-2)^i \,\, ], \end{aligned}%
\]
a quantity which is (essentially) equivalent to one arising in a conjecture of
Charney and Davis \cite{CharneyDavis}, related to a conjecture of Hopf (see
Section \ref{Hopf} below). Note that when $d$ is odd, $\sigma(P)$ vanishes by
the Dehn-Sommerville equations. When $d$ is even, we have the following result
(see Section \ref{signature}).

\begin{theorem}
\label{identify-signature} Let $P$ be a simple $d$-dimensional polytope, with
$d$ even. Then $\sigma(P)$ is the signature of the quadratic form $Q(x)=x^{2}$
defined on the $\frac{d}{2}^{th}$-graded component of McMullen's ring $\Pi(P)$.

In particular, when $P$ has rational vertices, $\sigma(P)$ is the signature or
index $\sigma(X_{\Delta})$ of the associated toric variety $X_{\Delta}$.
\end{theorem}

An important special case of the previously mentioned Charney-Davis conjecture
asserts that a certain combinatorial condition on $P$ (namely that of $\Delta$
being a \textit{flag} complex; see Section \ref{Hopf}) implies $(-1)^{d/2}%
\sigma(P) \geq0$. In this paper, we prove this conjecture when $P$ satisfies
certain stronger geometric conditions. We also give further conditions which
give lower bounds on $(-1)^{d/2}\sigma(P)$. In order to state these results,
we give rough definitions of some of these conditions here (see Section
\ref{Hirzebruch} for the actual definitions).

Say that the fan $\Delta$ is \textit{locally convex} (resp. \textit{locally
pointed convex, locally strongly convex}) if every $1$-dimensional cone in
$\Delta$ has the property that the union of all cones of $\Delta$ containing
it is convex (resp. pointed convex, strongly convex). For example (see
Propositions \ref{euclidean-imply-affine} and \ref{non-acute-by-polygons}
below), if each angle in every $2$-dimensional face of $P$ is non-acute (resp.
obtuse) then $\Delta$ will be locally convex (resp. locally strongly convex).
It turns out that $\Delta$ being locally convex implies that it is flag
(Proposition \ref{locally-convex-implies-flag}).

For a simple polytope $P$ with rational vertices, we define an integer $m(P)$
which measures how singular $X_{\Delta}$ is. To be precise, let $P$ in
$\mathbb{R}^{d}$ be rational with respect to some lattice $M$, and then $m(P)$
is defined to be the least common multiple over all $d$-dimensional cones
$\sigma$ of the normal fan $\Delta$ of the index $[N:N_{\sigma}]$, where $N$
is the lattice dual to $M$ and $N_{\sigma}$ is the sublattice spanned by the
lattice vectors on the extremal rays of $\sigma$. Note that the condition
$m(P)=1$ is equivalent to the smoothness of the toric variety $X_{\Delta}$,
and such polytopes $P$ are called \textit{Delzant} in the symplectic geometry
literature (e.g. \cite{Guillemin}).

Now we can state

\begin{theorem}
\label{lower-bounds} Let $P$ be a rational simple $d$-dimensional polytope
with $d$ even, and $\Delta$ its normal fan.

\begin{enumerate}
\item [(i)]If $\Delta$ is locally convex, then
\[
(-1)^{\frac{d}{2}}\sigma(P)\geq0.
\]

\item[(ii)] If $\Delta$ is locally pointed convex, then
\[
(-1)^{\frac{d}{2}}\sigma(P)\geq\frac{f_{d-1}(P)}{3m(P)^{d-1}}.
\]

\item[(iii)] If $\Delta$ is locally strongly convex, then
\[
(-1)^{\frac{d}{2}}\sigma(P)\geq\text{ coefficient of }x^{d}\text{ in }\left[
\frac{t^{d}}{m(P)^{d-1}}\,f(P,t^{-1})\right]  _{t\mapsto1-\frac{x}{\tan(x)}}.
\]
\end{enumerate}
\end{theorem}

\noindent We defer a discussion of the relation between Theorem
\ref{lower-bounds} (i) and the Charney-Davis conjecture to Section \ref{Hopf}.
It is amusing to see what Theorem \ref{lower-bounds} says beyond the
$g$-theorem, in the special case where $d=2$, that is, when $P$ is a
(rational) polygon. The $g$-theorem says exactly that
\[
f_{1}=f_{0} \geq3,
\]
or in other words, every polygon has the same number of edges as vertices, and
this number is at least $3$. Since
\[
(-1)^{\frac{d}{2}} \sigma(P) = f_{0}(P) - 4,
\]
Theorem \ref{lower-bounds} (i) tells us that when $\Delta$ is locally convex,
we must have $f_{0} \geq4$. In other words, triangles cannot have normal fan
$\Delta$ which is locally convex, as one can easily check. For $d=2$, the
conditions that $\Delta$ is locally pointed convex or locally strongly convex
coincide, and Theorem \ref{lower-bounds} (ii),(iii) both assert that under
these hypotheses, a (rational) polygon $P$ must have
\[
f_{0}(P) - 4 \,\, \geq\,\,\frac{f_{1}(P)}{3 m(P)} \,\,=\,\, \frac{f_{0}(P)}{3
m(P)}%
\]
or after a little algebra,
\begin{equation}
f_{0}(P) \,\, \geq\,\, \frac{12}{3-\frac{1}{m(P)}}.
\end{equation}
Since the right-hand side is strictly greater than $4$, we conclude that a
quadrilateral $P$ cannot have $\Delta$ locally pointed convex nor locally
strongly convex. This agrees with an easily-checked fact: a quadrilateral $P$
satisfies the weaker condition of having $\Delta$ locally convex if and only
if $P$ is a rectangle, and rectangles fail to have $\Delta$ locally pointed
convex. On the other hand, the inequality (\ref{polygon-inequality}) also
implies a not-quite-obvious fact: even though a (rational) pentagon can easily
have $\Delta$ locally strongly convex, this is impossible if $m(P)=1$, i.e.
there are no Delzant pentagons with this property. It is a fun exercise to
show directly that no such pentagon exists, and to construct a Delzant hexagon
with this property.

In fact, in the context of algebraic geometry, the proof of Theorem
\ref{lower-bounds} gives the following stronger assertion, valid for rational
simple polytopes of any dimension $d$ (not necessarily even) about the
expansion of the total $L$-class
\[
L(X) = L_{0}(X) + L_{1}(X) + \cdots+ L_{\frac{d}{2}}(X)
\]
where $L_{i}(X)$ is a cycle in $CH^{i}(X)_{\rationals}$, the Chow ring of $X$.

\begin{theorem}
\label{L-positive} Let $X=X_{\Delta}$ be a complete toric variety $X$
associated to a simplicial fan $\Delta$. If $\Delta$ is locally strongly
convex (resp. locally convex), then for each $i$ we have that $(-1)^{i}%
L_{i}(X)$ is effective (resp. either effective or $0$).
\end{theorem}

For instance, when $i=1$ this implies that if $\Delta$ is locally convex,
then
\[
\int_{X}(c_{1}^{2}(X)-2c_{2}(X))\cdot H_{1}\cdot\ldots\cdot H_{d-2}\leq0
\]
where $\{H_{i}\}$ are any ample divisor classes. This is reminiscent of the
Chern number inequality for the complex spinor bundle of $X$ when this bundle
is stable with respect to all polarizations; see e.g. \cite{Kobayashi}.

Notice that if $\Delta$ is not locally convex, $(-1)^{i}L_{i}(X)$ need not be
effective. For example, if $\Delta$ is the normal fan of the standard
$2$-dimensional simplex having vertices at $(0,0),(1,0),(0,1)$, then $X$ is
the complex projective plane, and $-L_{1}\left(  X\right)  $ is represented by
the negative of the Poincar\'{e} dual of a point.

\section{The alternating sum as signature}

\label{signature}

We wish to prove Theorem \ref{identify-signature}, whose statement we recall here.

\vskip.1in \noindent\textbf{Theorem \ref{identify-signature}.} \textit{Let $P$
be a simple $d$-dimensional polytope, with $d$ even. Then $\sigma(P)$ is the
signature of the quadratic form $Q(x)=x^{2}$, defined on the $\frac{d}{2}%
^{th}$-graded component of McMullen's ring $\Pi(P)$. }

\textit{In particular, when $P$ has rational vertices, $\sigma(P)$ is the
signature or index $\sigma(X_{\Delta})$ of the associated toric variety
$X_{\Delta}$. } \vskip.1in \begin{proof}
Taking $r=\frac{d}{2}$ in a result of McMullen
\cite[Theorem 8.6]{McMullen}, we find
that the quadratic form $(-1)^{\frac{d}{2}} Q(x)$ on the
$\frac{d}{2}^{th}$-graded component of $\Pi(P)$ has
$$
\begin{aligned}
\sum_{i=0}^{\frac{d}{2}} (-1)^i h_{\frac{d}{2}-i}(P)
&\text{ positive eigenvalues, and} \\
\sum_{i=0}^{\frac{d}{2}-1} (-1)^i h_{\frac{d}{2}-i-1}(P)
&\text{ negative eigenvalues.}
\end{aligned}
$$
Consequently, the signature $\sigma(Q)$ of $Q$ is
$$
\begin{aligned}
\sigma(Q)
&= (-1)^{\frac{d}{2}} \left[
\sum_{i=0}^{\frac{d}{2}} (-1)^i h_{\frac{d}{2}-i}(P)-
\sum_{i=0}^{\frac{d}{2}-1} (-1)^i h_{\frac{d}{2}-i-1}(P)
\right] \\
&= \sum_{i=0}^d (-1)^i h_i(P)
\end{aligned}
$$
where the second equality uses the {\it Dehn-Sommerville}
equations \cite[\S4]{McMullen}:
$$
h_i(P) = h_{d-i}(P).
$$
The second assertion of the theorem follows immediately
from McMullen's identification of the ring $\Pi(P)$ with
the quotient of the Stanley-Reisner ring of $\Delta$ by a certain
linear system of parameters \cite[\S14]{McMullen},
which is known to be isomorphic (after a doubling of the grading)
with the cohomology of $X_\Delta$ \cite[\S5.2]{Fulton}.
\end{proof}

\begin{remark}
\textrm{\ \newline Starting from any complete rational simplicial fan $\Delta
$, one can construct a toric variety $X_{\Delta}$ which will be complete, but
not necessarily projective, and satisfies Poincar\`{e} duality. The $h$-vector
for $\Delta$ can still be defined, and again has an interpretation as the
even-dimensional Betti numbers of $X_{\Delta}$ (see \cite[\S5.2]{Fulton}). We
suspect that the alternating sum of the $h$-vector is still the signature of
this complete toric variety. }

\textrm{Generalizing in a different direction, to \textit{any} polytope $P$
which is not necessarily simple, one can associate the normal fan $\Delta$ and
a projective toric variety $X_{\Delta}$. Although the (singular) cohomology of
$X_{\Delta}$ does not satisfy Poincar\'{e} duality, its \textit{intersection
cohomology} (in middle perversity) $IH^{\cdot}(X_{\Delta})$ will. There is a
combinatorially-defined \textit{generalized $h$-vector} which computes these
$IH^{\cdot}$ Betti numbers (see \cite{Stanley-generalized-h-vector}).
Moreover, using the Hard Lefschetz Theorem for intersection cohomology and the
fact that $X_{\Delta}$ is a finite union of affine subvarieties, the
alternating sum of the generalized $h$-vector equals the signature of the
quadratic form on $IH^{\cdot}(X_{\Delta})$ defined by the intersection
product. }
\end{remark}

\begin{remark}
\textrm{\ \newline The special case of the second assertion in Theorem
\ref{identify-signature} is known when $X_{\Delta}$ is smooth (i.e. $P$ is a
Delzant polytope); see \cite[Theorem 3.12 (3)]{Oda}. }
\end{remark}

\section{Lower bounds derived from the signature theorem}

\label{Hirzebruch}

The goal of this section is to explain the various notions used in Theorem
\ref{lower-bounds}, and to prove this theorem.

We begin with a $d$-dimensional lattice $M \cong\mathbb{Z}^{d}$ and its
associated real vector space $M_{\reals}=M\otimes_{\integers}\mathbb{R}$. A
\textit{polytope} $P$ in $M_{\reals}$ is the convex hull of a finite set of
points in $M_{\reals}$. We say that $P$ is \textit{rational} if these points
can be chosen to be rational with respect to the lattice $M$. The
\textit{dimension} of $P$ is the dimension of the smallest affine subspace
containing it. A \textit{face} of $P$ is the intersection of $P$ with one of
its supporting hyperplanes, and a face is always a polytope in its own right.
\textit{Vertices} and \textit{edges} of $P$ are $0$-dimensional and $1
$-dimensional faces, respectively. Every vertex of a $d$-dimensional polytope
lies on at least $d$ edges, and $P$ is called \textit{simple} if every vertex
lies on exactly $d$ edges.

Let $N=Hom(M,\mathbb{Z})$ be the dual lattice to $M$ and $N_{\reals}%
=N\otimes_{\integers}\mathbb{R}$ be the dual vector space to $M_{\reals}$,
with the natural pairing $M_{\reals}\otimes N_{\reals}\rightarrow\mathbb{R} $
denoted by $\langle\cdot, \cdot\rangle$. For a polytope $P$ in $M_{\reals}$,
the \textit{normal fan} $\Delta$ is the following collection of polyhedral
cones in $N_{\reals}$:
\[
\Delta= \{\sigma_{F}: F \text{ a face of }P\},
\]
where
\[
\sigma_{F}:=\{ v \in N_{\reals}: \langle u, v \rangle\leq\langle u^{\prime}, v
\rangle\text{ for all }u \in F, u^{\prime}\in P \}
\]
Note that

\begin{enumerate}
\item [$\bullet$]the normal fan $\Delta$ is a \textit{complete} fan, that is,
it exhausts $N_{\mathbb{R}}$,

\item[$\bullet$] $\Delta$ is a \textit{rational} fan, in the sense that its
rays all have rational slopes, if $P$ is rational,

\item[$\bullet$] if $P$ is $d$-dimensional, then every cone $\sigma_{F}$ in
$\Delta$ will be \textit{pointed}, that is, it will contain no proper
subspaces of $N_{\mathbb{R}}$,

\item[$\bullet$] $P$ is a simple polytope if and only if $\Delta$ is a
\textit{simplicial fan}, that is, every cone $\sigma$ in $\Delta$ is
simplicial in the sense that its extremal rays are linearly independent.
\end{enumerate}

We next define several affinely invariant conditions on a complete simplicial
fan $\Delta$ in $N_{\reals}$ (and hence on simple polytopes $P$ in $M_{\reals
}$) that appear in Theorem \ref{lower-bounds}. For any collection of
polyhedral cones $\Delta$ in $N_{\reals}$, let $|\Delta|$ denote the support
of $\Delta$, that is the union of all of its cones as a point set. Define the
\textit{star} and \textit{link} of one of the cones $\sigma$ in $\Delta$
similarly to the analogous notions in simplicial complexes: $\mathrm{star}%
_{\Delta}(\sigma)$ is the subfan consisting of those cones $\tau$ in $\Delta$
such that $\sigma, \tau$ lie in some common cone of $\Delta$, while
$\mathrm{link}_{\Delta}(\sigma)$ is the subfan of $\mathrm{star}_{\Delta
}(\sigma)$ consisting of those cones which intersect $\sigma$ only at the
origin. For a \textit{ray} (i.e. a $1$-dimensional cone) $\rho$ of $\Delta$,
say that the fan $\mathrm{star}_{\Delta}(\rho)$ is

\begin{enumerate}
\item [$\bullet$]\textit{convex} if its support $|\mathrm{star}_{\Delta}%
(\rho)|$ is a convex set in the usual sense,

\item[$\bullet$] \textit{pointed convex} if $|\mathrm{star}_{\Delta}(\rho)|$
is convex and contains no proper subspace of $N_{\mathbb{R}}$,

\item[$\bullet$] \textit{strongly convex} if furthermore for every cone
$\sigma$ in $\mathrm{link}_{\Delta}(\rho)$, there exists a linear hyperplane
$H$ in $N_{\mathbb{R}}$ which supports $\mathrm{star}_{\Delta}(\rho)$ and
whose intersection with $\mathrm{star}_{\Delta}(\rho)$ is exactly $\sigma$.
\end{enumerate}

Say that $\Delta$ is \textit{locally convex} (resp. \textit{locally pointed
convex}, \textit{locally strongly convex}) if every ray $\rho$ of $\Delta$ has
$\mathrm{star}_{\Delta}(\rho)$ convex (resp. pointed convex, strongly convex).
One has the easy implications
\[
\text{ locally strongly convex } \Rightarrow\text{ locally pointed convex }
\Rightarrow\text{ locally convex }.
\]

We recall here that the affine-lattice invariant $m(P)$ for a rational
polytope $P$ was defined (in the introduction) to be the least common multiple
of the positive integers $[N:N_{\sigma}]$ as $\sigma$ runs over all $d
$-dimensional cones in $\Delta$. Here $N_{\sigma}$ is the $d$-dimensional
sublattice of $N$ generated by the lattice vectors on the $d$ extremal rays of
$\sigma$. In a sense $m(P)$ measures how singular $X_{\Delta}$ is \cite[\S
2.6]{Fulton}, with $m(P)=1$ if and only if $X_{\Delta}$ is smooth, in which
case we say that $P$ is \textit{Delzant}.

We can now recall the statements of Theorems \ref{lower-bounds} and
\ref{L-positive}.

\vskip.1in \noindent\textbf{Theorem \ref{lower-bounds}.} \textit{Let $P$ be a
simple $d$-dimensional polytope in $M_{\reals}$, which is rational with
respect to $M$, and $\Delta$ its normal fan in $N_{\reals}$. Assume $d$ is
even. }

\begin{enumerate}
\item [(i)]\textit{If $\Delta$ is locally convex, then
\[
(-1)^{\frac{d}{2}}\sigma(P)\geq0.
\]
}

\item[(ii)] \textit{If $\Delta$ is locally pointed convex, then
\[
(-1)^{\frac{d}{2}}\sigma(P)\geq\frac{f_{d-1}(P)}{3m(P)^{d-1}}.
\]
}

\item[(iii)] \textit{If $\Delta$ is locally strongly convex, then
\[
(-1)^{\frac{d}{2}}\sigma(P)\geq\text{ coefficient of }x^{d}\text{ in }\left[
\frac{t^{d}}{m(P)^{d-1}}\,f(P,t^{-1})\right]  _{t\mapsto1-\frac{x}{\tan(x)}}.
\]
}
\end{enumerate}

\vskip.1in \noindent\textbf{Theorem \ref{L-positive}.} \textit{Let $X =
X_{\Delta}$ be a complete toric variety $X$ associated to a simplicial fan
$\Delta$, and let the expansion of the total $L$-class be
\[
L(X) = L_{0}(X) + L_{1}(X) + \cdots+ L_{\frac{d}{2}}(X)
\]
where $L_{i}(X)$ is a cycle in $CH^{i}(X)_{\rationals}$, the Chow ring of $X$. }

\textit{If $\Delta$ is locally strongly convex (resp. locally convex), then
for each $i$ we have that $(-1)^{i} L_{i}(X)$ is effective (resp. either
effective or $0$). } \vskip.1in

The remainder of this section is devoted to the proof of these theorems. We
begin by recalling some toric geometry. As a general reference for toric
varieties, we rely on Fulton \cite{Fulton}, although many of the facts we will
use can also be found in Oda's book \cite{Oda} or Danilov's survey article
\cite{Danilov}.

Let $X$ denote the toric variety $X_{\Delta}$. Simpleness of $P$ implies that
$X$ is an orbifold \cite[\S2.2]{Fulton}. Recall that irreducible toric
divisors\footnote{Actually these are $\mathbb{Q}$-Cartier divisors on the
orbifold $X$.} on $X$ correspond in a one-to-one fashion with the codimension
$1$ faces of $P$, or to $1$-dimensional rays in the normal fan $\Delta$.
Number these toric divisors on $X$ as $D_{1},...,D_{m}$. Intersection theory
for these $D_{i}$'s is studied in Chapter 5 of \cite{Fulton}. Every $D_{i}$ is
a toric variety in its own right with at worst orbifold singularities.
Moreover $D=\textstyle\bigcup_{i=1}^{m}D_{i}$ is a simple normal crossing
divisor on $X$ \cite[\S4.3]{Fulton}.

Next we want to express the signature of $X$ in terms of these $D_{i}$'s. When
$X$ is a smooth variety, a consequence of the hard Lefschetz Theorem is that
its signature $\sigma\left(  X\right)  $ can be expressed in terms of the
Hodge numbers of $X$ as follows \cite[Theorem 15.8.2]{Hirzebruch}:
\[
\sigma\left(  X\right)  = \Sigma_{p, q=0}^{d}\left(  -1\right)  ^{q}%
h^{p,q}\left(  X\right) .
\]
By the Dolbeault Theorem, $h^{p,q}\left(  X\right)  = \dim H^{q}(X,\Omega
_{X}^{p})$, and hence the signature can be expressed in terms of twisted
holomorphic Euler characteristics
\[
\sigma\left(  X\right)  =\sum_{p=0}^{d}\chi\left(  X,\Omega_{X}^{p}\right)
\]
where $\chi(X,E):=\sum_{q=0}^{d}(-1)^{q} \dim H^{q}(X,E)$. Using the
Riemann-Roch formula, we can write
\[
\chi\left(  X,\Omega_{X}^{p}\right)  =\int_{X}ch\left(  \Omega_{X}^{p}\right)
Td_{X}\text{,}%
\]
where $ch$ is the Chern character and $Td_{X}$ is the Todd class of $X$.
Therefore
\[
\sigma\left(  X\right)  =\int_{X}\textstyle\sum_{p=0}^{d}ch\left(  \Omega
_{X}^{p}\right)  Td_{X}.
\]

When $X$ is smooth, $\sum_{p=0}^{d}ch\left(  \Omega_{X}^{p}\right)  Td_{X}$
equals the Hirzebruch L-class $L\left(  X\right)  $ of $X$ (see page 16 in
\cite[Theorem 15.8.2]{Hirzebruch} for example) and we recover the Hirzebruch
Signature Formula
\begin{equation}
\sigma\left(  X\right)  =\int_{X}L\left(  X\right)  \text{.}\label{HSF}%
\end{equation}
If $X$ is a projective variety with at worst orbifold singularities, the hard
Lefschetz, Dolbeault, and Riemann-Roch Theorems continue to hold, and we can
take the sum $\sum_{p=0}^{d}ch\left(  \Omega_{X}^{p}\right)  Td_{X}$ as a
definition of $L\left(  X\right)  $. Since we can express $L\left(  X\right)
$ in terms of Chern roots of the orbi-bundle $\Omega_{X}^{1}$ and Chern
classes for orbi-bundles satisfy the same functorial properties as for the
usual Chern classes, the same holds true for $L\left(  X\right)  $. For
example we will use the splitting principle in the proof of the next lemma,
where we write $L\left(  X\right)  $ in terms of toric data.\footnote{For a
general orbifold $X$, not necessary an algebraic variety, Kawasaki
\cite{Kawasaki} expressed the signature of $X$ in terms of integral of certain
curvature forms, thus generalizing the Hirzebruch signature formula in a
different way.}

%Next we want to express the signature of $X$ in terms of these $D_{i}$'s.
%When $X$ is a smooth variety, its signature $\sigma \left( X\right) $ can be
%expressed in terms of the Hodge numbers of twisted holomorphic Euler characteristics, 
%\begin{equation*}
%\sigma \left( X\right) =\sum_{p=0}^{d}\chi \left( X,\Omega _{X}^{p}\right) .
%\end{equation*}
%Here $\chi \left( X,\Omega _{X}^{p}\right) :=
%\Sigma _{q=0}^{d}\left( -1\right) ^{q}h^{p,q}\left( X\right).$
%This formula follows from the hard
%Lefschetz theorem and therefore it continues to hold true for orbifolds.
%Using Riemann-Roch formula, we can write 
%\begin{equation*}
%\chi \left( X,\Omega _{X}^{p}\right) =\int_{X}ch\left( \Omega
%_{X}^{p}\right) Td_{X}\text{,}
%\end{equation*}
%and therefore 
%\begin{equation*}
%\sigma \left( X\right) =\int_{X}\textstyle\sum_{p=0}^{d}ch\left( \Omega
%_{X}^{p}\right) Td_{X}.
%\end{equation*}
%
%When $X$ is smooth, $\Sigma ch\left( \Omega _{X}^{p}\right) Td_{X}$ equals
%to the Hirzebruch L-class $L\left( X\right) $ of $X$ and we recover the
%Hirzebruch Signature Formula 
%\begin{equation}
%\label{HSF}
%\sigma \left( X\right) =\int_{X}L\left( X\right) \text{.}
%\end{equation}
%If $X$ is a projective variety with at worst orbifold singularities, we can
%take the sum $\Sigma ch\left( \Omega _{X}^{p}\right) Td_{X}$ as a definition of $%
%L\left( X\right) $. When we express $L\left( X\right) $ in terms of Chern
%roots of the orbi-bundle $\Omega _{X}$, they satisfy the same functorial
%properties for the usual L-class. 
%

\begin{lemma}
\label{crucial-calculation}
\[
\begin{aligned} (-1)^{\frac{d}{2}} \sigma\left( X\right) &= (-1)^{\frac{d}{2}}%
\int_{X}L\left( X\right) \\ & =\sum_{p=1}^{d/2} \,\, \sum_{\substack
{n_{1}+...+n_{p}=d/2\\n_{i}>0;i_{1}<...<i_{p}}} b_{n_{1}}\cdot\cdot\cdot
b_{n_{p}} \,\, \left(-1\right) ^{p}D_{i_{1}}^{2n_{1}} \cdot\ldots\cdot
D_{i_{p}}^{2n_{p}} \end{aligned}%
\]
where $D_{j_{1}}\cdot\ldots\cdot D_{j_{d}}$ denotes the intersection number
for the $d$ divisors $D_{j_{1}},\ldots,D_{j_{d}}$ on $X$, and $b_{n}$ are the
coefficients in the expansion
\[
\frac{\sqrt{x}}{\tanh\sqrt{x}}=1-\textstyle\sum_{n=1}^{\infty}\left(
-1\right)  ^{n}b_{n}x^{n}.
\]
That is, $b_{n}=\frac{2^{2n}B_{n}}{\left(  2n\right)  !}$ where $B_{n}$ is the
$n^{th}$ Bernoulli number.
\end{lemma}

\begin{proof}
Recall \cite{Milnor} that the L-class is a multiplicative
characteristic class corresponding to the power series $\frac{\sqrt{x}}%
{\tanh\sqrt{x}}$, as we are about to explain.
In our situation,
we need to compute the $L$-class of a holomorphic orbi-bundle $E$,
namely the tangent orbi-bundle of $X$.
For the purposes of this computation, we can treat
$E$ like a genuine vector bundle (see e.g.
\cite[Appendix A]{CoxKatz}).  We will assume that $E$ can be stably split into
a direct sum of line bundles, that is
$$
E \oplus O_X^{\oplus(m-d)} \cong\bigoplus_{i=1}^m L_i.
$$
Then $c(E) = \prod_{i=1}^m (1+x_i)$ where
$x_i$ is the first Chern class $c_1(L_i)$,
i.e. the $x_i$'s are {\it stable Chern roots} of $E$.
We then have
$$
c(E \otimes_{\reals} \complexes) = c(E \oplus\overline{E})
= \prod_{i=1}^m (1+x_i) \prod_{j=1}^m (1-x_j) = \prod_{i=1}^m (1- x_i^2).
$$
The $L$-class is then computed by the formula
$$
L(E) = \prod_{i=1}^m (1 - \sum_{n \geq1} (-1)^n b_n x_i^{2n})
$$
where $b_n$ is the positive number defined in the Lemma.
For example, in terms of Pontrjagin classes of $X$ we have
$$
\begin{aligned}
L_{1}\left(  X\right)  &=\frac{1}{3}p_{1}\\
L_{2}\left(  X\right) & =\frac{1}%
{45}\left(  7p_{2}-p_{1}^{2}\right)\\
L_{3}\left(  X\right)&=\frac{1}{945}\left( 62p_{3}-13p_{2}p_{1}+2p_{1}%
^{3}\right).
\end{aligned}
$$
To use the Hirzebruch Signature Formula \eqref{HSF},
we need to express $L\left(  X\right)  $ in
terms of toric data, specifically intersection numbers of the
toric divisors $D_{1},...,D_{m}$ on $X$ discussed above.
To relate these divisors with characteristic classes of $X$,
we consider the exact sequence of sheaves:%
\[
0\rightarrow\Omega_{X}^{1}\rightarrow\Omega_{X}^{1}\left(  \log D\right)
\overset{\text{Res}}{\longrightarrow}%
%TCIMACRO{\tbigoplus _{i=1}^{m}}%
%BeginExpansion
{\textstyle\bigoplus_{i=1}^{m}}
%EndExpansion
O_{D_{i}}\rightarrow0
\]
where $\Omega_{X}^{1}$ (resp. $\Omega_{X}^{1}\left(  \log D\right)  $) is the
sheaf of differentials on $X$
(resp. differentials on $X$ with logarithmic poles along $D$). Notice
that $\Omega_{X}^{1}\left(  \log D\right)  $ is a trivial sheaf of rank $d$
\cite[\S4.3]{Fulton}. On the other hand, there is an exact
sequence of sheaves%
\[
0\rightarrow O_{X}\left(  -D_{i}\right)  \rightarrow O_{X}\rightarrow
O_{D_{i}}\rightarrow0
\]
for each toric divisor $D_{i}$. From the functorial properties of Chern
classes, we have
$$
c\left(  \Omega_{X}^{1}\right)  =\Pi_{i=1}^{m}c\left(
O_{X}\left(  -D_{i}\right)  \right)  =\Pi_{i=1}^{m}\left(  1-D_{i}\right).
$$
Here we have identified a divisor $D_i$ with the Poincar\'e dual of the
first Chern class of its associated line bundle.  Since $\Omega_{X}^{1}$
is the sheaf of sections of the cotangent bundle on $X$,
we can write the total
Chern class of (the tangent bundle on) $X$ as
$$
c\left(  E\right)  =\Pi_{i=1}^{m}\left(  1+D_{i}\right).
$$
Namely these $D_{i}%
$'s behave as stable Chern roots of the tangent bundle of $X$.
Therefore, by the multiplicative property of the L-class, we have%
\[
L\left(  X\right)  =\prod_{i=1}^{m}\left(  1-%
%TCIMACRO{\textstyle\sum _{n=1}^{\infty}}%
%BeginExpansion
{\textstyle\sum_{n=1}^{\infty}}
%EndExpansion
\left(  -1\right)  ^{n}b_{n}D_{i}^{2n}\right)  .
\]
\noindent
Expanding the right-hand side of the above equality gives
the equality asserted in the lemma.
\end{proof}

To prove Theorem \ref{lower-bounds} we need to give lower bounds on the
intersection numbers
\[
\left(  -1\right)  ^{p}D_{i_{1}}^{2n_{1}}\cdot\cdot\cdot D_{i_{p}}^{2n_{p}}%
\]
that appear in the right-hand-side of Lemma \ref{crucial-calculation}, under
our various hypotheses on the fan $\Delta$. We begin by rewriting
\begin{equation}
\left(  -1\right)  ^{p}D_{i_{1}}^{2n_{1}}\cdot\cdot\cdot D_{i_{p}}^{2n_{p}%
}=\int_{D_{i_{1}}\cap...\cap D_{i_{p}}}\left(  -D_{i_{1}}\right)  ^{2n_{1}%
-1}\cdot\cdot\cdot\left(  -D_{i_{p}}\right)  ^{2n_{p}-1}.
\end{equation}
This expression leads us to consider the restriction of the line bundles
$O_{X}(-D_{i})$ to the subvarieties $D_{i_{1}}\cap...\cap D_{i_{p}}$. For any
of the irreducible toric divisors $D_{i}$ on $X$, let $O_{D_{i}}(-D_{i})$
denote the restriction of $O_{X}(-D_{i})$ to the toric subvariety $D_{i}$
(this is the \textit{conormal bundle} of $D_{i}$ in $X$). Recall that for an
invertible sheaf $O(E)$ on a $d$-dimensional orbifold $X$, one says that
$O(E)$ is \textit{big} if the corresponding divisor $E$ satisfies $E^{d} > 0$.
The key observation in obtaining the desired lower bounds is then

\begin{lemma}
\label{ample-big-nef} Let $P$ be a simple $d$-dimensional polytope in
$M_{\mathbb{R}}$, which is rational with respect to $M$, and $\Delta$ its
normal fan in $N_{\mathbb{R}}$. Let $D_{i}$ be any of the irreducible toric
divisors on $X=X_{\Delta}$.

\begin{enumerate}
\item [(i)]If $\Delta$ is locally convex, then $O_{D_{i}}(-D_{i})$ is
generated by global sections.

\item[(ii)] If $\Delta$ is locally pointed convex, then $O_{D_{i}}(-D_{i})$ is
generated by global sections and big.

\item[(iii)] If $\Delta$ is locally strongly convex, then $O_{D_{i}}(-D_{i})$
is ample.
\end{enumerate}
\end{lemma}

Assuming Lemma \ref{ample-big-nef} for the moment, we finish the proof of
Theorem \ref{lower-bounds}.

\vskip.1in \noindent\textit{Proof of Theorem \ref{lower-bounds} (i)}. Under
the assumption that $\Delta$ is locally convex, we know that the restriction
of $O_{X}\left(  -D_{i_{j}}\right)  $ to $D_{i_{1}}\cap...\cap D_{i_{p}}$ is
generated by global sections for $1\leq j\leq p$ by Lemma \ref{ample-big-nef}.
This implies that the integral (\ref{integral}) equals the intersection number
of such divisors on the toric subvariety $D_{i_{1}}\cap...\cap D_{i_{p}}$ and
therefore it is nonnegative \footnote{This follows from the fact that the
divisor class of a line bundle which is generated by global sections is a
limit of $\mathbb{Q}$-divisors which are ample. Positivity of intersection
numbers of ample divisors is well-known; see e.g. \cite[Chapter 12]%
{Fulton-intersection-theory}} , that is
\[
\left(  -1\right)  ^{p}D_{i_{1}}^{2n_{1}}\cdot\cdot\cdot D_{i_{p}}^{2n_{p}%
}\geq0.
\]
The non-negativity asserted in Theorem \ref{lower-bounds} (i) now follows
term-by-term from the sum in Lemma \ref{crucial-calculation}.$\qed$

\vskip.1in \noindent\textit{Proof of Theorem \ref{lower-bounds} (ii)}. If
$\Delta$ is locally pointed convex, then $O_{D_{i}}\left(  -D_{i}\right)  $ is
generated by global sections and big. The bigness of $O_{D_{i}}\left(
-D_{i}\right)  $ on $D_{i}$ implies that
\[
-D_{i}^{d}=\int_{D_{i}}\left(  -D_{i}\right)  ^{d-1}%
\]
is strictly positive.

\noindent\textbf{Claim.} $-D_{i}^{d} \geq\frac{1}{m(P)^{d-1}}$.

\noindent To prove this, we proceed as in the \textit{algebraic moving lemma}
\cite[\S5.2, p. 107]{Fulton}, making repeated use of the fact that if $n_{j} $
is the first non-zero lattice point on the ray of $\Delta$ corresponding to
$D_{j}$, then for any $u$ in $M$, one has
\begin{equation}
\sum_{j}\langle u,n_{j}\rangle D_{j}=0\label{moving-equation}%
\end{equation}
in the Chow ring \cite[Proposition, Part (ii), \S5.2, p. 106]{Fulton}. This
allows one to take intersection monomials that contain some divisor $D_{j_{0}%
}$ raised to a power greater than $1$, and replace one factor of $D_{j_{0}}$
by a sum of other divisors. By repeating this process for all of the monomials
in a total of $d-1$ stages, one can replace $D_{i}^{d}$ by a sum of the form
$\sum a_{j_{1},\ldots,j_{d}}D_{j_{1}}\cdots D_{j_{d}}$ in which each term has
$D_{j_{1}},\ldots,D_{j_{d}}$ \textit{distinct} divisors which intersect at an
isolated point of $X$, and each $a_{j_{1},\ldots,j_{d}}$ is a rational number.
We must keep careful track of the denominators of the coefficients introduced
at each stage.

At the first stage, by choosing any $u$ in $M$ with $\langle u,n_{i}\rangle
=1$, we can use \eqref{moving-equation} to replace one factor of $D_{i}$ in
$D_{i}^{d}$ by a sum of other divisors $D_{j}$ with \textit{integer
coefficients} (that is, introducing \textit{no} denominators). However, in
each of the next $d-2$ stages, when one wishes to use \eqref{moving-equation}
to substitute for a divisor $D_{j_{0}}$, one must choose $u$ in $M$
constrained to vanish on normal vectors $n_{j}$ for other divisors $D_{j}$ in
the monomial, and this may force the coefficient $\langle u,n_{j_{0}}\rangle$
of $D_{j_{0}}$ to be larger than $1$ in \eqref{moving-equation}, although it
will always be an integer factor of $m(P)$. Consequently, at each stage after
the first, we may introduce factors into the denominators that divide into
$m(P)$. Since there are $d-2$ stages after the first, we conclude that each
$a_{j_{1},\ldots,j_{d}}$ can be written with the denominator $m(P)^{d-2}$.
Finally, each intersection product $D_{j_{1}}\cdots D_{j_{d}}$ is the
reciprocal of the multiplicity at the corresponding point of $X$, which is the
index $[N:N_{\sigma}]$ where $\sigma$ is the $d$-dimensional cone of $\Delta$
corresponding to that point \cite[\S2.6]{Fulton}. Since each $[N:N_{\sigma}]$
divides $m(P)$, we conclude that $-D_{i}^{d}$ lies in $\frac{1}{m(P)^{d-1}%
}\mathbb{Z}$, and since it is positive, it is at least $\frac{1}{m(P)^{d-1}}$.

We have shown then that each term with $p=1$ on the right-hand side of Lemma
\ref{crucial-calculation} is at least $\frac{1}{m(P)^{d-1}}$, and the number
of such terms is the number of codimension one faces of $P$, i.e. $f_{d-1}(P)
$. Moreover we still have nonnegativity of the other terms $\left(  -1\right)
^{p}D_{i_{1}}^{2n_{1}}\cdot\cdot\cdot D_{i_{p}}^{2n_{p}}$ because $O_{D_{i}%
}\left(  -D_{i}\right)  $ is generated by global sections. Therefore, since
$b_{1}=\frac{1}{3}$, we conclude from Lemma \ref{crucial-calculation} that%

\[
\left(  -1\right)  ^{d/2}\sigma\left(  \Delta\right)  \geq\frac{f_{d-1}}{3
m(P)^{d-1}}. \qed
\]

\vskip.1in \noindent\textit{Proof of Theorem \ref{lower-bounds} (iii)}. If
$\Delta$ is locally strongly convex, then $O_{D_{i}}\left(  -D_{i}\right)  $
is ample. By similar arguments as in assertions (i) and (ii), we have
\[
\left(  -1\right)  ^{p}D_{i_{1}}^{2n_{1}}\cdot\cdot\cdot D_{i_{p}}^{2n_{p}%
}\geq\frac{1}{m(P)^{d-1}}%
\]
provided that $D_{i_{1}}\cap...\cap D_{i_{p}}$ is non-empty. By the simplicity
of $P$, each of its codimension $p$ faces can be expressed uniquely as the
intersection of distinct codimension one faces, corresponding to the non-empty
intersection of divisors $D_{i_{1}}, \ldots, D_{i_{p}}.$ Therefore, after
choosing positive integers $n_{1},\ldots, n_{p}$, the number of non-vanishing
terms of the form $\left(  -1\right)  ^{p}D_{i_{1}}^{2n_{1}}\cdot\cdot\cdot
D_{i_{p}}^{2n_{p}}$ in the expansion of Lemma \ref{crucial-calculation} is
$f_{d-p}(P)$. Hence%

\[
\begin{aligned} (-1)^{\frac{d}{2}} \sigma\left( X\right) &= (-1)^{\frac{d}{2}}%
\int_{X}L\left( X\right) \\ & =\sum_{p=1}^{d/2}\left( -1\right) ^{p}%
\sum_{\substack{n_{1}+...+n_{p}=d/2\\n_{i}>0;i_{1}<...<i_{p} }}b_{n_{1}}%
\cdot\cdot\cdot b_{n_{p}}D_{i_{1}}^{2n_{1}}\cdot\cdot\cdot D_{i_{p}}^{2n_{p}%
}\\ & \geq\sum_{p=1}^{d/2} \sum_{\substack{n_{1}+...+n_{p}=d/2\\n_{i}>0}%
}b_{n_{1}}\cdot\cdot\cdot b_{n_{p}} \frac{f_{d-p}(P)}{m(P)^{d-1}}%
\\ & = \sum_{p=1}^{d/2} \frac{f_{d-p}(P)}{m(P)^{d-1}} \left[ \,\,\text
{ coefficient of }x^d\text{ in } \left( \sum_{n \geq1} b_n x^{2n}%
\right)^p \,\, \right]. \end{aligned}%
\]

Note that
\[
\frac{\sqrt x}{\tanh\sqrt x} = 1 - \sum_{n \geq1} (-1)^{n} b_{n} x^{n}%
\]
implies
\[
\sum_{n \geq1} b_{n} x^{2n} = 1-\frac{x}{\tan(x)},
\]
and note also that
\[
\sum_{p \geq1} f_{d-p}(P) \, t^{p} = t^{d} f(P,t^{-1}).
\]
This allows us to rewrite the above inequality as in the assertion of Theorem
\ref{lower-bounds} (iii). $\qed$

\vskip.1in \noindent\textit{Proof of Theorem \ref{L-positive} }. Recall from
the proof of Lemma \ref{crucial-calculation} that the total $L$-class has
expansion
\[
L(X) = \sum_{p \geq1} \left(  -1\right)  ^{p}\sum_{\substack{(n_{1}%
,\ldots,n_{p}) \\n_{i}>0;i_{1}<...<i_{p}}} (-1)^{\sum n_{i}}b_{n_{1}}%
\cdot\cdot\cdot b_{n_{p}}D_{i_{1}}^{2n_{1}}\cdot\cdot\cdot D_{i_{p}}^{2n_{p}}.
\]
Consequently,
\[
(-1)^{j} L_{j}(X) = \sum_{p \geq1} \left( -1\right)  ^{p} \sum
_{\substack{n_{1}+\cdots+n_{p} = j \\n_{i}>0;i_{1}<...<i_{p}}} b_{n_{1}}%
\cdot\cdot\cdot b_{n_{p}}D_{i_{1}}^{2n_{1}}\cdot\cdot\cdot D_{i_{p}}^{2n_{p}}.
\]
Therefore it suffices to show that each term
\[
\left( -1\right)  ^{p}D_{i_{1}}^{2n_{1}} \cdot\cdot\cdot D_{i_{p}}^{2n_{p}}%
\]
is effective if $\Delta$ is locally strongly convex (the case where $\Delta$
is locally convex is similar). Here we use the fact from Lemma
\ref{ample-big-nef} that restriction of $O(-D_{i_{k}})$ to $D_{i_{k}}$ is
ample, and therefore also ample when further restricted to the transverse
intersection $V=D_{i_{1}} \cap\cdots\cap D_{i_{p}}$. Consequently, the cycle
class
\[
(-D_{i_{1}})^{2n_{1}-1} \cdot\cdot\cdot(-D_{i_{p}})^{2n_{p}-1}%
\]
is effective in the Chow ring $CH(V)_{\rationals}$ by Bertini's Theorem.
Therefore $\left( -1\right)  ^{p}D_{i_{1}}^{2n_{1}} \cdot\cdot\cdot D_{i_{p}%
}^{2n_{p}}$ is effective in $CH(X)_{\rationals}$. $\qed$

\vskip.1in \noindent\textit{Proof of Lemma \ref{ample-big-nef} }. We recall
some facts about toric divisors contained generally in \cite[\S3.3,
3.4]{Fulton}. In general, any divisor $E$ on $X$ defines a continuous
piecewise linear function $\Psi_{E}^{X}$ on the support $|\Delta|= N_{\reals}%
$. Every divisor $E$ on $X$ is linearly equivalent to a linear combination of
irreducible toric divisors. If we write $E=\textstyle\sum_{i=1}^{m}a_{i}D_{i}%
$, then $\Psi_{E}^{X}$ is determined by the property that $\Psi_{E}^{X}\left(
n_{i}\right)  =-a_{i}$ where $n_{i}$ is the first nonzero lattice point of $N$
along $\rho_{i}$. In particular, $\Psi_{-D_{i}}^{X}:N_{\reals}\rightarrow
\mathbb{R}$ is determined by $\Psi_{-D_{i}}^{X}\left(  n_{j}\right)
=\delta_{ij}$. The ampleness of the line bundle $O_{X}\left(  E\right)  $ can
be measured by the convexity of the piecewise linear function $\Psi_{E}^{X}$.
More precisely, $O_{X}\left(  E\right)  $ is ample (resp. generated by global
sections) if and only if $\Psi_{E}^{X}$ is strictly convex (resp. convex).

We now discuss assertions (i), (iii) of the lemma, leaving (ii) for later.
First we examine the particular case of the discussion in the previous
paragraph where the bundle is $O_{D_{i}}(-D_{i})$ on the toric subvariety
$D_{i}$. Assume the divisor $D_{i}$ corresponds to a ray $\rho_{i}$ in the
normal fan $\Delta$. The fan $\Delta^{D_{i}}$ associated to $D_{i}$ naturally
lives in the quotient space $N_{\reals}/ \rho_{i}$ (here we are abusing
notation by letting $\rho_{i}$ denote both a ray and also the $1$-dimensional
subspace spanned by this ray) \cite[\S3.1]{Fulton}. Then every cone in
$\Delta^{D_{i}}$ corresponds to a cone in $\Delta$ containing $\rho_{i}$ as a
face (and vice-versa), that is, a cone in $\mathrm{star}_{\Delta}(\rho_{i})$.
The boundary of $\mathrm{star}_{\Delta}(\rho_{i})$ is $\mathrm{link}_{\Delta
}(\rho_{i})$, and here we use the fact, proven in the Appendix, that
$\mathrm{link}_{\Delta}(\rho_{i})$ is affinely equivalent to the graph of the
continuous piecewise linear function $\Psi_{-D_{i}}^{D_{i}}: N_{\reals}/
\rho_{i} \rightarrow\mathbb{R}$. From the discussion in the previous
paragraph, we conclude that $O_{D_{i}}(-D_{i})$ is generated by global
sections (resp. ample) if $\Delta$ is locally convex (resp. locally strongly convex).

Lastly we discuss asertion (ii) of the lemma. We want to prove that $O_{D_{i}%
}\left(  -D_{i}\right) $ is big for every irreducible toric divisor $D_{i}$
under the assumption that $\Delta$ is locally pointed convex. The fact that
$\Delta$ is locally pointed convex says that the space $|\mathrm{star}%
_{\Delta}(\rho_{i})|$ is a pointed convex polyhedral cone. There is then a
unique fan $\Sigma$ having the following properties:

\begin{enumerate}
\item [$\bullet$]$\Sigma$ is refined by $\mathrm{star}_{\Delta}(\rho_{i})$,
and they have the same support, that is,
\[
|\Sigma|=|\mathrm{star}_{\Delta}(\rho_{i})|,
\]

\item[$\bullet$] $\rho_{i}$ is the only ray in the interior of $\Sigma$, and

\item[$\bullet$] $\Sigma$ is strongly convex in the sense that that every ray
of $\Sigma$ except for $\rho_{i}$ is the intersection of $|\Sigma|$ with some
supporting hyperplane.
\end{enumerate}

This cone $\Sigma$ projects to a complete fan $\bar{\Delta}^{D_{i}}$ in
$N_{\reals}/ \rho_{i}$, which is refined by $\Delta^{D_{i}}$. Therefore we
obtain a birational morphism \cite[\S1.4]{Fulton}.
\[
\pi:D_{i}=X_{ \Delta^{D_{i}}} \rightarrow X_{\bar{\Delta}^{D_{i}}}.
\]
Moreover $\Psi_{-D_{i}}^{D_{i}}$ still defines a continuous piecewise linear
function on $N_{\reals}/ \rho_{i}$ , the support of $\bar{\Delta}^{D_{i}}$,
which is now strongly convex. Therefore it defines an ample Cartier divisor on
$X_{\bar{\Delta}^{D_{i}}} $. Call this divisor $C$. Then it is not difficult
to see that $O_{D_{i}}\left(  -D_{i}\right)  $ is just the pullback of
$O_{X_{\bar{\Delta}^{D_{i}}} }\left(  C\right)  $. Moreover, their
self-intersection numbers are equal, that is
\[
\left(  -D_{i}\right)  ^{d-1}=C^{d-1}\text{.}%
\]
Now $C$ is an ample divisor on $X_{\bar{\Delta}^{D_{i}}} $ and therefore
$C^{d-1}$ is strictly positive. Hence the same is true for $-D_{i}$. That is,
$O_{D_{i}}(-D_{i})$ is a big line bundle on the toric subvariety $D_{i}$. This
completes the proof of Lemma \ref{ample-big-nef}. $\qed$

The previous proof raises the following question: is the assumption of
rationality for the simple polytope $P$ really necessary in Theorem
\ref{lower-bounds}? In approaching this problem, it would be interesting if
Lemma \ref{crucial-calculation} and the intersection numbers that appear
within it have some interpretation purely within the convexity framework used
by McMullen \cite{McMullen}.

\section{Examples}

\label{examples}

In this section we discuss examples of simple polytopes $P$ whose normal fans
$\Delta$ satisfy the hypotheses of Theorem \ref{lower-bounds}.

We begin with some properties of $P$ that are Euclidean invariants, so we will
assume that $M_{\reals}$ is endowed with a (positive definite) inner product
$\langle\cdot, \cdot\rangle$ which identifies $M_{\reals}$ with its dual space
$N_{\reals}$. Thus we can think of both $P$ and its normal fan $\Delta$ as
living in $M_{\reals}$.

Say that a polytope $P$ is \textit{non-acute in codimension $1$} (resp.
\textit{obtuse in codimension $1$}) if every codimension $2$ face of $P$ has
the property that the dihedral angle between the two codimension $1$ faces
meeting there is non-acute (resp. obtuse), that is, at least (resp. greater
than) $\frac{\pi}{2}$. Say that $P$ is \textit{non-acute} (resp.
\textit{obtuse}) if $P$ and every one of its faces of each dimension
considered as polytopes in their own right are non-acute (resp. obtuse) in
codimension $1$. We have the following obvious implications:%
\[%
\begin{array}
[c]{ccc}%
& \text{obtuse} & \\
\Longrightarrow & \text{obtuse in codimension }1\text{ and non-acute} & \\
\Longrightarrow & \text{non-acute} &
\end{array}
\]

The next proposition shows that these Euclideanly invariant conditions on $P$
imply the affinely invariant conditions on $\Delta$ defined in the Section
\ref{Hirzebruch}.

\begin{proposition}
\label{euclidean-imply-affine} Let $P$ be a simple $d$-dimensional polytope in
$M_{\reals}$, with normal fan $\Delta$.

\begin{enumerate}
\item [(i)]If $P$ is non-acute, then $\Delta$ is locally convex.

\item[(ii)] If $P$ is obtuse in codimension $1$ and non-acute, then $\Delta$
is locally pointed convex.

\item[(iii)] If $P$ is obtuse, then $\Delta$ is locally strongly convex.
\end{enumerate}
\end{proposition}

The next corollary then follows immediately from the previous proposition and
Theorem \ref{lower-bounds}.

\begin{corollary}

\label{euclidean-imply-bounds} Let $P$ be a simple, rational $d$-dimensional
polytope in $M_{\reals}$, with normal fan $\Delta$.

\begin{enumerate}
\item [(i)]If $P$ is non-acute, then
\[
(-1)^{\frac{d}{2}} \sigma(P) \geq0.
\]

\item[(ii)] If $P$ is obtuse in codimension $1$, and non-acute, then
\[
(-1)^{\frac{d}{2}} \sigma(P) \geq\frac{f_{d-1}(P)}{3m(P)^{d-1}}.
\]

\item[(iii)] If $P$ is obtuse, then
\[
(-1)^{\frac{d}{2}} \sigma(P) \geq\text{ coefficient of }x^{d}\text{ in }
\left[  \frac{t^{d}}{m(P)^{d-1}} \, f(P,t^{-1}) \right] _{t \mapsto1-\frac
{x}{\tan(x)}}.
\]
\end{enumerate}
\end{corollary}

\begin{remark}
\textrm{\ \newline It is easy to see that obtuse simple polytopes can always
be made rational without changing their facial structure by a slight
perturbation of their facets, so that one might think of removing the
rationality assumption from part (iii) of the previous corollary. However,
after this perturbation it is not clear what the lattice-invariant $m(P)$ will
be, i.e. it could be any positive integer. }

\textrm{It is not obvious whether a non-acute, simple polytope always has the
same facial structure as a rational, non-acute, simple polytope. This would
follow if every non-acute, simple polytope had the same facial structure as an
obtuse, simple polytope, but this is false. For example, a regular
$3$-dimensional cube is non-acute and simple, but no obtuse polytope can have
the facial structure of a $3$-cube. }
\end{remark}

\begin{remark}
\textrm{\ \newline \label{non-acute-bound-without-rationality} \noindent M.
Davis has pointed out to us that the first assertion of Corollary
\ref{euclidean-imply-bounds} can be proven using facts from
\cite{CharneyDavis} and the mirror construction $M(P)$ of the next section,
without any assumption that $P$ is rational. We defer a sketch of this proof
until the description of $M(P)$ at the end of that section. }
\end{remark}

\vskip.1in \noindent\textit{Proof of Proposition \ref{euclidean-imply-affine}%
.} We begin by rephrasing some of our definitions about non-acuteness and
obtuseness in terms of $\Delta$. Obtuseness (resp. non-acuteness) in
codimension $1$ for $P$ corresponds to the following property of $\Delta$: any
two vectors $n, n^{\prime}$ spanning the extremal rays of a $2$-dimensional
cone of $\Delta$ must have
\[
\langle n,n^{\prime}\rangle> 0 \qquad( \text{ resp. } \langle n,n^{\prime
}\rangle\geq0 ).
\]
Similarly, obtuseness (resp. non-acuteness) for $P$ corresponds to the
following property of $\Delta$: any vectors $n_{1}, \ldots, n_{t}$ spanning
the extremal rays of a $t$-dimensional (simplicial) cone of $\Delta$ must
have
\[
\langle\pi(n_{1}),\pi(n_{2}) \rangle> 0 \qquad( \text{ resp. } \langle
\pi(n_{1}),\pi(n_{2}) \rangle\geq0 )
\]
where $\pi$ is the orthogonal projection onto the space perpendicular to the
span of the vectors $n_{3}, n_{4}, \ldots, n_{t}$.

Having said this, observe that if $P$ is non-acute in codimension $1$, for any
ray $\rho$ in $\Delta$ (spanned by a vector which we name $n$), the hyperplane
$\rho^{\perp}$ normal to $\rho$ supports $\mathrm{star}_{\Delta}(\rho) $: if
$P$ is non-acute in codimension $1$, we must have $\langle n,n^{\prime}%
\rangle\geq0$ for each vector $n^{\prime}$ spanning a ray in $\mathrm{star}%
_{\Delta}(\rho)$, and hence for every vector in $\mathrm{star}_{\Delta}(\rho
)$. Similarly, if $P$ is obtuse in codimension $1$ then this hyperplane
$\rho^{\perp}$ not only supports $\mathrm{star}_{\Delta}(\rho)$, but also
intersects it only in the origin. Consequently, assertion (ii) of the lemma
follows once we prove assertion (i).

For assertions (i), (iii), we make use of the fact that strong or weak
convexity of $\mathrm{star}_{\Delta}(\rho)$ can be checked \textit{locally} in
a certain way, similar to checking regularity of triangulations (see e.g.
\cite[\S1.3]{DeLoera}). Roughly speaking, each cone $\sigma$ in the link of
$\rho$ must have the property that the union of cones containing $\sigma$
within $\mathrm{link}_{\Delta}(n)$ ``bend outwards" at $\sigma$ away from
$\rho$, rather than ``bending inward" toward $\rho$. To be more formal,
consider every minimal dependence of the form
\begin{equation}
\sum_{i \in F} \alpha_{i} n_{i} = \beta n + \sum_{j \in G} \beta_{j} m_{j}%
\end{equation}
where

\begin{enumerate}
\item [-]$\{n_{i}\}_{i\in F}$ are vectors spanning the extremal rays of some
cone $\sigma$ in $\mathrm{link}_{\Delta}(n)$,

\item[-] each $m_{j}$ for $f$ in $G$ spans a ray in $\mathrm{link}_{\Delta
}(\sigma)$,

\item[-] the coefficients $\alpha_{i}, \beta_{j}$ are all strictly positive.
\end{enumerate}

Then $\mathrm{star}_{\Delta}(\rho)$ is strictly convex if and only in every
such dependence we have $\beta< 0$. It is weakly convex if and only if in
every such dependence we have $\beta\leq0$.

As a step toward proving assertions (i), (iii), given a dependence as in
(\ref{minimal-dependence}) we apply the orthogonal projection $\pi$ onto the
space perpendicular to all of the $\{n_{i}\}_{i \in F}$, yielding the
following equation
\[
0 = \beta\, \pi(n) + \sum_{j \in G} \beta_{j} \, \pi(m_{j}),
\]
and then taking the inner product with $\pi(n)$ on both sides yields
\begin{equation}
0 = \beta\, \langle\pi(n), \pi(n) \rangle+ \sum_{j \in G} \beta_{j} \,
\langle\pi(m_{j}) , \pi(n) \rangle.
\end{equation}

To prove (iii), we assume $P$ is obtuse and that there is some choice of a
dependence as in (\ref{minimal-dependence}) such that $\beta\geq0$. But then
we reach a contradiction in Equation (\ref{contradiction-equation}), because
we assumed $\beta_{j} > 0$, we have $\langle\pi(m_{j}) , \pi(n) \rangle> 0$ by
virtue of the obtuseness of $P$, and $\langle\pi(n), \pi(n) \rangle$ is always non-negative.

To prove (i), we assume $P$ is non-acute and that there is some choice of a
dependence as in (\ref{minimal-dependence}) such that $\beta> 0$. Then similar
considerations in equation (\ref{contradiction-equation}) imply that we must
have $\langle\pi(n), \pi(n) \rangle= 0$, i.e. $\pi(n) = 0$. This would imply
$\langle n_{i}, n \rangle=0$ for each $i$ in $F$. To reach a contradiction
from this, take the inner product with $n$ on both sides of equation
(\ref{minimal-dependence}), to obtain
\[
0= \beta\, \langle n, n \rangle+ \sum_{j \in G} \beta_{j} \, \langle m_{j}, n
\rangle.
\]
Non-acuteness (even in codimension 1) of $P$ implies $\langle m_{j}, n
\rangle\geq0$, and $\langle n, n \rangle$ is always positive, so this last
equation is a contradiction to $\beta> 0$. $\qed$

One source of non-acute simple polytopes are finite Coxeter groups (see
\cite[Chapter 1]{Humphreys} for background). Recall that a finite
\textit{Coxeter group} is a finite group $W$ acting on a Euclidean space and
generated by reflections. Given a finite Coxeter group $W$, there is
associated a set of (normalized) \textit{roots} $\Phi$ by taking all the unit
normals of reflecting hyperplanes. Let $Z$ be the \textit{zonotope} (
\cite[\S7.3]{Ziegler}) associated with $\Phi$, that is,
\[
Z= \left\{  \sum_{\alpha\in\Phi} c_{\alpha}\alpha: 0 \leq c_{\alpha}%
\leq1\right\} .
\]

\begin{proposition}
\label{Coxeter-zonotopes-non-acute} The zonotope $Z$ associated to any finite
Coxeter group $W$ is non-acute and simple. Furthermore $Z$ is obtuse in
codimension $1$ if $W$ is irreducible.
\end{proposition}

\begin{proof}
We refer to \cite{Humphreys} for all facts about Coxeter groups
used in this proof.
By general facts about zonotopes \cite[\S7.3]{Ziegler},
the normal fan $\Delta$ of $Z$ is the complete fan cut out
by the hyperplanes associated with reflections in $W$.
The maximal cones in this fan are the {\it Weyl chambers} of $W$, which
are all simplicial cones.  Hence $Z$ is a simple polytope.
To show that $Z$ is non-acute, we must show that each of its
faces is non-acute in codimension $1$.  However, these faces are
always affine translations of Coxeter zonotopes corresponding to
standard parabolic subgroups of $W$.  So we only need to show
$Z$ itself is non-acute in codimension $1$.  This
is equivalent to showing that every pair of rays in $\Delta$
which span a $2$-dimensional cone form a non-obtuse angle.
Because $W$ acts transitively on the
Weyl chambers in $\Delta$, we may assume that this
pair of rays lie in the {\it fundamental Weyl chamber}, that is,
we may assume that these rays come from the dual basis
to some choice of simple roots $\alpha_1,\ldots,\alpha_d$.
Since every choice of simple roots has the property that
$\langle\alpha_i, \alpha_j \rangle\leq0$ for all $i \neq j$,
the first assertion follows from the first part of
Lemma \ref{pairwise-non-acute} below.
The second assertion follows from Lemma \ref{pairwise-non-acute} (ii)
below.  This is because the obtuseness graph for any choice of
simple roots associated with a Coxeter group $W$
is isomorphic to the (unlabelled) Coxeter graph,
and the Coxeter graph is connected exactly when $W$ is irreducible.
\end{proof} 

The following lemma was used in the preceding proof.

\begin{lemma}
\label{pairwise-non-acute} Let $\{\alpha_{i}\}_{i=1}^{d}$ be a basis for
$\mathbb{R}^{d}$ with $\langle\alpha_{i}, \alpha_{j} \rangle\leq0$ for all $i
\neq j$. Then the dual basis $\{\alpha^{\vee}_{i}\}_{i=1}^{d}$ satisfies

\begin{enumerate}
\item [(i)]$\langle\alpha^{\vee}_{i}, \alpha^{\vee}_{j} \rangle\geq0$ for all
$i \neq j$, and

\item[(ii)] $\langle\alpha^{\vee}_{i}, \alpha^{\vee}_{j} \rangle> 0$ for all
$i \neq j$ if the ``obtuseness graph" on $\{1,2,\ldots,d\}$, having an edge
$\{i,j\}$ whenever $\langle\alpha_{i}, \alpha_{j} \rangle< 0$, is connected.
\end{enumerate}
\end{lemma}

\begin{proof}
We prove assertion (i) by induction
on $d$, with the cases $d=1,2$ being trivial.
In the inductive step, assume $d \geq3$.  Without loss of
generality, we must show $\langle\alpha_1, \alpha_2 \rangle\geq0$.
Let $\pi: \reals^d \rightarrow\alpha_d^\perp$ be
orthogonal projection.  Write
$$
\alpha_i = \pi(\alpha_i) + c_i \alpha_d
$$
for each $i \leq d-1$.
Our first claim is that $c_i \leq0$ for each $i \leq d-1$.
To see this, note that
$$
\begin{aligned}
0 &\geq\langle\alpha_i, \alpha_d \rangle\\
&= \langle\pi(\alpha_i), \alpha_d \rangle+
c_i \langle\alpha_d, \alpha_d \rangle\\
&= c_i \langle\alpha_d, \alpha_d \rangle.
\end{aligned}
$$
Our second claim is that
$\langle\pi(\alpha_i), \pi(\alpha_j) \rangle\leq0$
for $1 \leq i \neq j \leq d-1$.
To see this, note that
$$
\begin{aligned}
0 &\geq\langle\alpha_i, \alpha_j \rangle\\
&= \langle\pi(\alpha_i), \pi(\alpha_j) \rangle+
c_j \langle\alpha_d, \pi(\alpha_j) \rangle+
c_i \langle\pi(\alpha_i), \alpha_d \rangle+
c_i c_j \langle\alpha_d, \alpha_d \rangle\\
&=\langle\pi(\alpha_i), \pi(\alpha_j) \rangle+
c_i c_j \langle\alpha_d, \alpha_d \rangle.
\end{aligned}
$$
and the last term in the last sum is non-negative by our
first claim.
Our third claim is that $\{\pi(\alpha_i)\}_{i=1}^{d-1}$
and  $\{\alpha^\vee_i\}_{i=1}^{d-1}$ are dual bases inside
$\alpha_d^\perp$.  To see this, note that
$$
\begin{aligned}
\delta_{ij} &= \langle\alpha_i, \alpha^\vee_j \rangle\\
&=\langle\pi(\alpha_i), \alpha^\vee_j \rangle+
c_i \langle\alpha_d, \alpha^\vee_j \rangle\\
&= \langle\pi(\alpha_i), \alpha^\vee_j \rangle.
\end{aligned}
$$
From the second and third claims, we can apply induction
to conclude that $\langle\alpha^\vee_i, \alpha^\vee_j \rangle\geq0$
for $1 \leq i \neq j \leq d-1$, and in particular this holds
for $i=1,j=2$ as desired.
To prove assertion (ii), we use the same induction on $d$, with
the cases $d=1,2$ still being trivial.  We must in addition show that
if $\{\alpha_i\}_{i=1}^d$ have connected obtuseness graph, then
there is a re-indexing (that is a choice of $\alpha_d$)
so that $\{\pi(\alpha_i)\}_{i=1}^{d-1}$ also satisfies this hypothesis.
To achieve this, let $\alpha_d$ correspond to a node $d$ in the
obtuseness graph whose removal does not disconnect it, e.g. choose
$d$ to be a leaf in some spanning tree for the graph.
Then for $i \neq j$ with $i,j \leq d-1$ we have
$$
\langle\pi(\alpha_i) , \pi(\alpha_j) \rangle=
\langle\alpha_i , \alpha_j \rangle-
c_i c_j \langle\alpha_d , \alpha_d \rangle.
$$
This implies $\pi(\alpha_i) , \pi(\alpha_j)$ were obtuse whenever
$\alpha_i, \alpha_j$ were, so the obtuseness graph remains connected.
\end{proof}

\begin{remark}
\textrm{\ \newline If the finite Coxeter group $W$ is
\textit{crystallographic} (or a \textit{Weyl} group) then a crystallographic
root system associated with $W$ gives a more natural choice of hyperplane
normals to use than the unit normals in defining the Coxeter zonotope $Z$.
With this choice, the normal fan $\Delta$ is not only rational with respect to
the \textit{weight lattice} $N$, but also $m(Z)=1$ with respect to the dual
lattice $M$. Hence $Z$ is Delzant, so that the toric variety $X_{\Delta}$ is
smooth. }
\end{remark}

For the classical Weyl groups $W$ of types $A,B(=C),D$, there are known
generating functions for the $h$-vectors of the associated Coxeter zonotopes
$Z$, which specialize to give explicit generating functions for the signature
$\sigma(Z)$. The $h$-vector in this case turns out to give the distribution of
the elements of the Weyl group $W$ according to their \textit{descents}, i.e.
the number of simple roots which they send to negative roots (see
\cite{Bjorner-Coxetercomplexes}). Generating functions for the descent
distribution of all classical Weyl groups may be found in \cite{Reiner}. For
example, it follows from these that if $Z_{A_{n-1}}$ is the Coxeter zonotope
of type $A_{n-1}$, then we have the formula
\[
\sum_{n \geq0} \sigma(Z_{A_{n-1}}) \frac{x^{n}}{n!} = \tanh(x)
\]
which was computed in \cite[Example p. 52]{EdelmanReiner} for somewhat
different reasons.

The fact that Coxeter zonotopes have locally convex normal fans also follows
because these normal fans come from \textit{simplicial hyperplane
arrangements} (we thank M. Davis for suggesting this). Say that an arrangement
of linear hyperplanes ${}$ in $\mathbb{R}^{d}$ is \textit{simplicial} if it
decomposes $\mathbb{R}^{d}$ into a simplicial fan.

\begin{proposition}
\label{simple-zonotope-implies-locally-convex} The fan $\Delta$ associated to
a simplicial hyperplane arrangement ${}$ is locally convex.
\end{proposition}

\begin{proof}
For each ray $\rho$ of $\Delta$, we will express $\str_\Delta(\rho)$ as
an intersection of closed half-spaces defined by a subset of the
hyperplanes of ${\mathcal A}$, thereby showing that it is convex.
To describe this intersection, note that since $\Delta$ is simplicial, given
any chamber ($d$-dimensional cone) $\sigma$ of
$\Delta$ that contains $\rho$, there is a unique hyperplane
$H_\sigma$ bounding $\sigma$ which does not contain $\rho$.  Choose
a linear functional $u_\sigma$ which vanishes on $H_\sigma$
and is positive on $\rho$, and then we claim that
$$
|\str_\Delta(\rho)| =
\bigcap_{\text{chambers }\sigma\supset\rho} \{u_\sigma\geq0\}.
$$
To see that the left-hand side is contained in the right,
note that for any chamber $\sigma$ containing $\rho$ and any hyperplane
$H$ in ${\mathcal A}$ not containing $\rho$, we must have $\sigma$ and
$\rho$ on the same side of $H$.  Consequently, for every pair of
chambers $\sigma, \sigma'$ containing $\rho$ we have $u_\sigma\geq0$
on $\sigma'$ (and symmetrically $u_{\sigma'} \geq0$ on $\sigma$).
This implies the desired inclusion.
To show that the right-hand side is contained in the left, since both
sets are closed and $d$-dimensional, it suffices to show that every
chamber in the left is contained in the right, or contrapositively, that
every chamber not contained in the right is not in the left.  Given a
chamber $\sigma$ not in the right, consider the unique chamber $\sigma'$
containing $\rho$ which is ``perturbed in the direction of $\sigma
$''.  In other
words, $\sigma'$ is chosen so that it
contains a vector $v+\epsilon w$ where $v$ is any non-zero vector
in $\rho$, $w$ is any vector in the interior of $\sigma$, and $\epsilon$ is a
very small positive number.  Since $\sigma$ does not contain $\rho$, we know
$\sigma\neq\sigma'$, and hence there is at least one hyperplane of ${\mathcal
A}$
separating them.  Since $\Delta$ is simplicial, every bounding hyperplane
of $\sigma'$ except for $H_{\sigma'}$ will contain $r$, and hence have
$\sigma$ and $\sigma'$ on the same side (by construction of $\sigma'$).  This
means $H_{\sigma'}$ must separate $\sigma$ and $\sigma'$, so $u_{\sigma'} < 0$
on $\sigma$, implying $\sigma$ is not in the left-hand side.
\end{proof}

The Coxeter zonotopes of type $A$ are related to another infinite family of
simple polytopes, the \textit{associahedra}, which turn out to have locally
convex normal fans. Recall \cite{Lee} that the associahedron $_{n}$ is an
$(n-3)$-dimensional polytope whose vertices correspond to all possible
parenthesizations of a product $a_{1} a_{2} \cdots a_{n-1}$, and having an
edge between two parenthesizations if they differ by a single ``rebracketing".
Equivalently, vertices of $_{n}$ correspond to triangulations of a convex
$n$-gon, and there is an edge between two triangulations if they differ only
by a "diagonal flip" within a single quadrilateral.

\begin{proposition}
The associahedron $_{n}$ has a realization as a simple convex polytope whose
normal fan $\Delta_{n}$ is locally convex.
\end{proposition}

\begin{proof}
In \cite[\S3]{Lee}, the normal fan $\Delta_n$ is thought of as a simplicial
complex, and more precisely, as the
boundary of a simplicial polytope $Q_n$ having the
origin in its interior.  There $Q_n$ is
constructed by a sequence of stellar subdivisions
of certain faces of an $(n-3)$-simplex having vertices labelled
$1,2,\ldots,n-2$.   Since the normal fan $\Delta_n$ is simplicial,
the associahedron is simple (as is well-known).
Our strategy for showing $\Delta_n$ is locally convex
is to relate it to the Weyl chambers of
type $A_{n-3}$. If we assume that the $(n-3)$-simplex above is regular, and
take its barycenter as the origin in $\reals^{n-3}$, then the barycentric
subdivision of its boundary is a simplicial polytope isomorphic to
the Coxeter complex for type $A_{n-3}$.  Hence the normal fan
$\Delta_n$ of $\assoc_n$ refines the fan of Weyl chambers for type $A_{n-3}$.
Note that an alternate description of this Weyl chamber fan is that it is the
set of all chambers cut out by the hyperplanes $x_i = x_j$, that is,
its (open) chambers are defined by inequalities of the form
$x_{\pi_1} > x_{\pi_2} >\cdots>x_{\pi_{n-2}}$
for permutations $\pi$ of $\{1,2,\ldots,n-2\}$.
To show $\Delta_n$ is locally convex, we must first identify the rays
$\rho$ of $\Delta_n$, and then show that $\str_{\Delta_n}(\rho)$ is a pointed
convex cone.  According to the construction of \cite[\S3]{Lee},
a ray $\rho_{ij}$ of $\Delta_n$ corresponds to the barycenter of
a face of the $(n-3)$-simplex which
is spanned by a set of vertices labelled by a {\it contiguous}
sequence $i,i+1,\ldots,j-1,j$ with $1 \leq i \leq j \leq n-2$,
with $(i,j) \neq(1,n-2)$.
It is then not hard to check from the construction that
$\str_\Delta(\rho_{ij})$ consists of the union of all (closed) chambers
for type $A_{n-3}$ which satisfy the inequalities
$$
x_{i},x_{i+1},\ldots,x_{j-1},x_{j} \geq x_{i-1}, x_{j+1}
$$
(where here we omit the inequalities involving $x_{i-1}$ if $i=1$,
and similarly for $x_{j+1}$ if $j=n+2$).  It is clear that these
inequalities describe a convex cone, and hence $\Delta_n$ is
locally convex.
\end{proof}

\noindent It follows then from this Proposition and Theorem \ref{lower-bounds}%
(ii) that $(-1)^{\frac{n-3}{2}} \sigma(_{n}) \geq0$ for $n$ odd (and of
course, $\sigma(_{n}) = 0$ for $n$ even). However, as in the case of Coxeter
zonotopes of type $A$, we can compute $\sigma(_{n})$ explicitly using the
formulas for the $f$-vector or $h$-vector of $_{n}$ given in \cite[Theorem
3]{Lee}. Specifically, these formulas imply that for $n \geq3$ we have
\[
\begin{aligned} \sigma(\mathcal A_n) &=\sum_{i = 0}^{n-3} (-1)^i \frac{1}{n-1}%
\binom{n-3}{i} \binom{n-1}{i+1} \\ &={}_2F_1 \left( \left. \begin{matrix}%
3-n & 2-n \\ & 2 \end{matrix} \right| -1 \right)\\ &=\left\{ \begin{matrix}%
(-1)^{\frac{n-3}{2}}C_{\frac{n-1}{2}} & \text{ if }n\text{ is odd}%
\\ 0 & \text{ if }n\text{ is even} \end{matrix} \right. \end{aligned}%
\]
where $C_{n}$ denotes the Catalan number $\frac{1}{n}\binom{2n-2}{n-1}$. Here
the ${}_{2}F_{1}$ is hypergeometric series notation, and the last equality
uses Kummer's summation of a well-poised $_{2}F_{1}$ at $-1$ (see e.g.
\cite[p. 9]{Bailey}).

Returning to the discusion of non-acute and obtuse polytopes, it is worth
noting the following facts, pointed out to us by M. Davis. Recall that a
simplicial complex $K$ is called \textit{flag} if every set of vertices
$v_{1},\ldots, v_{r}$ which pairwise span edges of $K$ also jointly span an
$(r-1) $-simplex of $K$.

\begin{proposition}
\label{non-acute-by-polygons} A polytope $P$ is non-acute (resp. obtuse) if
and only if each of its $2$-dimensional faces are non-acute (resp. obtuse).

Furthermore, non-acuteness of any polytope $P$ implies that $P$ is simple
\end{proposition}

\begin{proof}
The first assertion for non-acute polytopes
follows from an easy lemma due to Moussong \cite[Lemma 2.4.1]{CharneyDavis}.
In the notation of \cite{CharneyDavis},
saying that  every $2$-dimensional face of
$P$ is non-acute (in codimension $1$) is equivalent to saying that
for every vertex $v$ of $P$, the spherical
simplex $\sigma= Lk(v,P)$ has size $\geq\frac{\pi}{2}$.  Then
\cite[Lemma 2.4.1]{CharneyDavis}
asserts that every link $Lk(\tau,\sigma)$ of a face of this spherical
simplex also has size $\geq\frac{\pi}{2}$.  But a face $F$ of
$P$ containing $v$ has $Lk(v,F)$ of the form $Lk(\tau,\sigma)$
for some $\tau$, and hence $F$ is non-acute in codimension $1$
when considered as a polytope in its own right.  That is, $P$ is non-acute.
An easy adaptation of this argument to the obtuse case
proves the first assertion of the proposition for obtuse polytopes.
The fact that non-acuteness implies simplicity again comes from considering
the spherical simplex $\sigma= Lk(v,P)$ for any vertex $v$, which will
have size $ \geq\frac{\pi}{2}$.  Then its polar dual
spherical convex polytope $\sigma^*$ will have all of its dihedral
angles less than or equal to $\frac{\pi}{2}$.  This forces $\sigma^*$
to be a spherical simplex, by \cite[p. 44]{Vinberg}, and hence $\sigma$
itself must be a spherical simplex.  This implies $v$ has exactly $d$ neighbors,
so $P$ is simple.
\end{proof}

Obtuse polytopes turn out to be relatively scarce in comparison with non-acute
polytopes, For example, it is easily seen that Coxeter zonotopes, although
always non-acute by Proposition \ref{Coxeter-zonotopes-non-acute}, are not in
general obtuse in dimensions $3$ and higher. It is easy to find obtuse
polytopes in dimensions up to $4$:

\begin{enumerate}
\item [$\bullet$]in dimension 2, the regular $n$-gons for $n \geq5$,

\item[$\bullet$] in dimension 3, the dodecahedron,

\item[$\bullet$] in dimension 4 the ``$120$-cell" (see \cite[pp.
292-293]{Coxeter})
\end{enumerate}

However M. Davis has pointed out to us that in dimensions $5$ higher, there
are no obtuse polytopes, due to a result of Kalai \cite[Theorem 1]{Kalai} (see
also \cite[p. 68]{Vinberg} for the case of simple polytopes): every
$d$-dimensional convex polytope for $d \geq5$ contains either a triangular or
quadrangular $2$-dimensional face.

\section{Relation to conjectures of Hopf and of Charney and Davis}

\label{Hopf}

In this section we discuss the relation of Theorem \ref{lower-bounds}(i) to
the conjectures of Hopf and of Charney and Davis mentioned in the Introduction.

Let $M^{d}$ be a compact $d$-dimensional closed Riemannian manifold. When $d$
is odd, Poincar\'e duality implies that the Euler characteristic $\chi(M)$
vanishes. When $d$ is even, a conjecture of H. Hopf (see e.g.
\cite{CharneyDavis}) asserts that if $M^{d}$ has non-positive sectional
curvature, the Euler characteristic $\chi(M^{d})$ satsfies
\[
(-1)^{\frac{d}{2}} \chi(M^{d}) \geq0.
\]
This result is known for $d=2,4$ by Chern's Gauss-Bonnet formula, but open for
general $d$; see \cite[\S0]{CharneyDavis} for some history.

Charney and Davis \cite{CharneyDavis} explored a combinatorial analogue of
this conjecture, and we refer the reader to their paper for terms which are
not defined precisely here. Let $M^{d}$ be a compact $d$-dimensional closed
manifold which has the structure of a (locally finite) \textit{Euclidean cell
complex}, that is, it is formed by gluing together convex polytopes via
isometries of their faces. One can endow such a cell complex with a metric
space structure that is Euclidean within each polytopal cell, making it a
\textit{geodesic space}. Gromov has defined a notion of when a geodesic space
is \textit{nonpositively curved}, and Charney and Davis made the following conjecture:

\begin{conjecture}
\cite[Conjecture A]{CharneyDavis} \label{conjecture-A} If $M^{d}$ is a
non-positively curved, piecewise Euclidean, closed manifold with $d$ even,
then
\[
(-1)^{\frac{d}{2}}\chi(M^{d})\geq0.
\]
\end{conjecture}

For piecewise Euclidean cell complexes, nonpositive curvature turns out to be
equivalent to a local condition at each vertex. Specifically, at each vertex
$v$ of $M^{d}$, one has a piecewise spherical cell complex $Lk(v,M^{d})$
called the \textit{link} of $v$ in $M^{d}$, which is homeomorphic to a
\textit{generalized homology $(d-1)$-sphere}, and inherits its own geodesic
space structure. Nonpositive curvature of $M^{d}$ turns out to be a metric
condition on each of these complexes $Lk(v,M^{d})$. Charney and Davis show
\cite[(3.4.3)]{CharneyDavis} that the Euler characteristic $\chi(M^{d})$ can
be written as the sum of certain local quantities $\kappa( Lk(v,M^{d}) )$
defined in terms of the metric structure on $Lk(v,M^{d})$:%

\begin{equation}
\chi(M^{d}) = \sum_{v} \kappa( Lk(v,M^{d}) ).
\end{equation}

In the special case where the polytopes in the Euclidean cell decomposition of
$M^{d}$ are all right-angled cubes, the links $Lk(v,M^{d})$ are all simplicial
complexes, and the quantity $\kappa( Lk(v,M^{d}) )$ has a simple combinatorial
expression purely in terms of the numbers of simplices of each dimension in
these complexes (that is, independent of their metric structure). Furthermore,
in this case, non-positive curvature corresponds to the combinatorial
condition that each link is a \textit{flag complex}, that is, the minimal
subsets of vertices in $Lk(v,M^{d})$ which do not span a simplex always have
cardinality two. They then noted that in this special case, Conjecture
\ref{conjecture-A} would follow via Equation (\ref{combinatorial-Gauss-Bonnet}%
) from

\begin{conjecture}
\cite[Conjecture D]{CharneyDavis} \label{conjecture-D} If $\Delta$ is a flag
simplicial complex triangulating a generalized homology $(d-1)$-sphere with
$d$ even, then
\[
(-1)^{\frac{d}{2}}\kappa(\Delta)\geq0.
\]
\end{conjecture}

\noindent This \textit{Charney-Davis conjecture} is trivial for $d=2$, has
recently been proven by Davis and Okun \cite{DavisOkun} for $d=4$ using
$L_{2}$-homology of Coxeter groups, and is also known (by an observation of
Babson and a result of Stanley; see \cite[\S7]{CharneyDavis}) for the special
class of flag simplicial complexes which are barycentric subdivisions of
boundaries of convex polytopes.

Local convexity of simplicial fans turns out to be stronger than flagness:

\begin{proposition}
\label{locally-convex-implies-flag} A locally convex complete simplicial fan
$\Delta$ in $\mathbb{R}^{d}$ is flag, when considered as a simplicial complex
triangulating a $(d-1)$-sphere.
\end{proposition}

\begin{proof}
Assume that $\Delta$ is not flag, so that there exist rays
$\rho_1,\ldots,\rho_k$ whose convex hull $\sigma:=\conv(\rho_1,\ldots,\rho_k)$
is {\it not} a cone of $\Delta$, but $\conv(\rho_i,\rho_j)$ {\it is} a cone
of $\Delta$ for each $i,j$.  Choose such a collection of rays
of minimum cardinality $k$, so that $\conv(\rho_1,\ldots,\hat{\rho_i}%
,\ldots,\rho_k)$
is a cone of $\Delta$ for each $i$ (in other words, the boundary complex
$\partial\sigma$ is a subcomplex of $\Delta$).
We wish to show that $\str_\Delta(\rho
_1)$ is not convex.  To see this, consider
$\sigma\cap\str_\Delta(\rho_1)$, that is, the collection of cones
$$
\{ \sigma' \cap\sigma: \sigma' \in\str_\Delta(\rho_1)\}.
$$
Since $\sigma$ is not a cone of $\Delta$ but $|\Delta|=\reals^d$,
this collection must contain at least one $2$-dimensional cone
of the form $\sigma' \cap\sigma= \conv(\rho_1,\rho)$,
where $\rho$ is a ray of $\sigma$ but $\rho\not\in\{\rho_2,\ldots,\rho_k\}$.
Since $\rho$ lies inside $\sigma$ and $\partial\sigma$ is a
subcomplex of $\Delta$, $\rho$ cannot lie in $\partial\sigma$ (else
some cone of $\partial\Delta$ would be further
subdivided, and not be a cone of $\Delta$).  Consequently $\rho$ lies
in the interior of $\sigma$.  Then the ray $\rho': =\rho- \epsilon\rho_1$
for very small $\epsilon> 0$ has the following properties:
\begin{enumerate}
\item[$\bullet$] $\rho'$ lies in $\sigma$, because $\rho
$ was in the interior of
$\sigma$,
\item[$\bullet$] $\rho'$ therefore lies in the convex hull of $\str
_\Delta(\rho_1)$,
since $\sigma$ does (as its extreme rays $\rho_1,\ldots,\rho_k$ of $\sigma$
are all in $\str_\Delta(\rho_1)$),
\item[$\bullet$] $\rho'$ does not lie in $\str_\Delta(\rho
_1)$, else it would lie
in a cone $\sigma'$ of $\Delta$ containing $\rho_1$, and then $\sigma'$ would
contain $\rho$ in the relative interior of one of its faces, a contradiction.
\end{enumerate}
Therefore $\str_\Delta(\rho_1)$ is not convex.
\end{proof}

In light of the preceeding proposition, one might ask if every flag simplicial
sphere has a realization as a locally convex complete simplicial fan. We thank
X. Dong for the following argument showing that an even weaker statement is
false. One can show that complete simplicial fans always give rise to
$PL$-spheres. Therefore if one takes the barycentric subdivision of any
regular cellular sphere which is not $PL$ (such as the double suspension of
Poincar\'{e}'s famous homology sphere), this will give a flag simplicial
sphere which is not $PL$ and therefore has no realization as a complete
simplicial fan (let alone one which is locally convex).

Our results were motivated by the Charney-Davis conjecture and the following
fact: when $P$ is a simple $d$-polytope and $\Delta$ is its normal fan
considered as a $(d-1)$-dimensional simplicial complex, one can check that
\begin{equation}
\sigma(P) = 2^{d} \kappa( \Delta).
\end{equation}

As a consequence, we deduce the following from Proposition
\ref{locally-convex-implies-flag} and Theorem \ref{lower-bounds} (i).

\begin{corollary}
\label{non-acute-lower-bound} Let $P$ be rational simple polytope, $\Delta$
its normal fan. If $\Delta$ is locally convex, then it is flag and satisfies
the Charney-Davis conjecture.

In particular by Corollary \ref{euclidean-imply-bounds}, if $P$ is a non-acute
simple rational polytope then its normal fan $\Delta$ is flag and satisfies
the Charney-Davis conjecture.
\end{corollary}

It is worth mentioning that the special case of Conjecture \ref{conjecture-A}
considered in \cite{CharneyDavis} where $M^{d}$ is decomposed into
right-angled cubes is ``polar dual" to another special case that fits nicely
with our results. Say that $M^{d}$ has a \textit{corner decomposition} if the
local structure at every vertex in the decomposition is combinatorially
isomorphic to the coordinate orthants in $\mathbb{R}^{d}$, that is, each link
$Lk(v,M^{d})$ has the combinatorial structure of the boundary complex of a
$d$-dimensional \textit{cross-polytope} or \textit{hyperoctahedron}. (Note
that this condition immediately implies that each of the $d$-dimensional
polytopes in the decomposition must be simple). A straightforward counting
argument (essentially equivalent to the calculation proving \cite[(3.5.2)]%
{CharneyDavis}) shows that for a manifold $M^{d}$ with corner decomposition
into simple polytopes $P_{1},\ldots,P_{N}$ one has
\begin{equation}
\chi(M^{d}) = \frac{1}{2^{d}} \sum_{i=1}^{N} \sigma(P_{i}).
\end{equation}

The following corollary is then immediate from this relation and Theorem
\ref{lower-bounds}.

\begin{corollary}
\label{corner-decomposition-lower-bound} Let $M^{d}$ be an $d$-dimensional
manifold with $d$ even, having a corner decomposition.

If each of the simple $d$-polytopes in the corner decomposition is rational
and has normal fan which is locally convex, then
\[
(-1)^{\frac{d}{2}}\chi(M^{d})\geq0.
\]
In particular, this holds if each of the simple $d$-polytopes is non-acute.
\end{corollary}

Several interesting examples of manifolds with corner decompositions into
simple polytopes that are either Coxeter zonotopes (hence non-acute) or
associahedra (hence locally convex) may be found in
\cite{DavisJanuszkiewiczScott}.

There is also an important general construction of such manifolds called
\textit{mirroring} which we now discuss. This construction (or its polar dual)
appears repeatedly in the work of Davis \cite{Davis-Annals, Davis-Duke,
DavisJanuszkiewicz, DavisJanuszkiewiczScott}, and was used in \cite[\S
6]{CharneyDavis} to show that the case of their Conjecture \ref{conjecture-A}
for manifolds decomposed into right-angled cubes is equivalent to their
Conjecture \ref{conjecture-D}. In a special case, this construction begins
with a generalized homology $(d-1)$-sphere $L$ with $n$ vertices and produces
a cubical orientable generalized homology $d$-manifold $ML$ having $2^{n}$
vertices, with the link at each of these vertices isomorphic to $L$. Hence we
have
\[
\chi(ML) = 2^{n} \cdot\kappa(L).
\]
We wish to make use of the polar dual of this construction, which applies to
an arbitrary simple $d$-dimensional polytope $P$, yielding an orientable
$d$-manifold $M(P)$ with a corner decomposition having every $d$-dimensional
cell isometric to $P$. The construction is as follows: denote the
$(d-1)$-dimensional faces of $P$ by $F_{1},F_{2},\ldots,F_{n}$, and let $M(P)$
be the quotient of $2^{n}$ disjoint copies $\{P_{\epsilon}\}_{\epsilon
\in\{+,-\}^{n}}$ of $P$, in which two copies $P_{\epsilon}, P_{\epsilon
^{\prime}}$ are identified along their face $F_{i}$ whenever $\epsilon,
\epsilon^{\prime}$ differ in the $i^{th}$ coordinate and nowhere else. As a
consequence of equation (\ref{local-formula}) we have
\begin{equation}
\chi(M(P)) = 2^{n-d} \cdot\sigma(P),
\end{equation}
which shows that the ``non-acute" assertions in Corollaries
\ref{corner-decomposition-lower-bound} and \ref{non-acute-lower-bound} are equivalent.

We can now use the mirror construction to complete the proof of an assertion
from the previous section. We are indebted to M. Davis for the statement and
proof of this assertion.

\vskip.1in \noindent\textit{Proof of Corollary \ref{euclidean-imply-bounds}
(i) without assuming rationality of $P$ (as referred to in Remark
\ref{non-acute-bound-without-rationality})}: Assume that $P$ is a simple
non-acute $d$-dimensional polytope with $d$ even. We wish to show that
$(-1)^{\frac{d}{2}} \sigma(P) \geq0$.

Construct $M(P)$ as above, a manifold with corner decomposition into non-acute
simple polytopes having $\chi(M(P)) = 2^{n-d} \sigma(P)$ if $P$ had $n$
codimension $1$ faces. In the notation of \cite{CharneyDavis}, this means that
all the links $Lk(v,M(P))$ have size $\geq\frac{\pi}{2}$ and are
combinatorially isomorphic to boundaries of cross-polytopes. This implies that
these links' underlying simplicial complexes are flag complexes satisfying
\cite[Conjecture D']{CharneyDavis}, and then \cite[Proposition 5.7]%
{CharneyDavis} implies that each of these links $Lk(v,M(P))$ satisfies
\cite[Conjecture C']{CharneyDavis}. This implies that $(-1)^{\frac{d}{2}}
\kappa( Lk(v,M(P)) ) \geq0$. Combining this with Equation
(\ref{combinatorial-Gauss-Bonnet}), we conclude that $(-1)^{\frac{d}{2}}
\chi(M(P)) \geq0$, and finally via Equation (\ref{mirror-equation}), that
$(-1)^{\frac{d}{2}} \sigma(P) \geq0$.

We note that a similar argument (involving an adaptation of \cite[Lemma
2.4.1]{CharneyDavis}) proves $(-1)^{\frac{d}{2}} \sigma(P) > 0$ when $P$ is
obtuse, but does not yield in any obvious way the stronger assertion of
Corollary \ref{euclidean-imply-bounds} (iii). $\qed$

\section{Appendix: the conormal bundle of a toric divisor}

In this appendix we describe the conormal bundle of a toric divisor on
$X=X_{\Delta}$ when the fan $\Delta$ is complete and simplicial. Denote the
collection of toric divisors on $X$ by $D_{1},\ldots,D_{m}.$ As mentioned in
Section \ref{Hirzebruch}, the conormal bundle of a divisor, say $D_{1}$, can
be identified as the restriction of $O_{X}(-D_{1})$ to $D_{1}$, which we
renamed $O_{D_{1}}\left(  -D_{1}\right) $. It corresponds to a continuous
piecewise linear function:
\[
\Psi_{-D_{1}}^{D_{1}}:N_{\mathbb{R}}/\rho_{1}\rightarrow\mathbb{R}%
\]
as in the discussion of Section \ref{Hirzebruch}. Here $\rho_{1}$ is the ray
in the fan $\Delta$ corresponding to the divisor $D_{1}$. We wish to identify
the graph of $\Psi_{-D_{1}}^{D_{1}}$ with $\mathrm{link}_{\Delta}\left(
\rho_{1}\right)  $, which we recall is the boundary of $\mathrm{star}_{\Delta
}(\rho_{1})$, the latter being the union of all cones of $\Delta$ containing
$\rho_{1}$.

\begin{proposition}
Let $D_{1}$ be a toric divisor of a toric variety $X=X_{\Delta}$ with $\Delta$
simplicial. Then the graph of the piecewise linear function for $O_{D_{1}%
}(D_{1}) $ is affinely equivalent to the boundary $link_{\Delta}\left(
\rho_{1}\right)  $ of $star_{\Delta}\left(  \rho_{1} \right) $, where
$\rho_{1}$ is the ray corresponding to $D_{1}$.
\end{proposition}

\begin{proof} We can index the toric divisors $D_1,\ldots,D_m$
of $X$ in such a way
that $D=D_{1}$ and $D_{2},...,D_{l}$ are those which are adjacent to $D_{1}$%
. Let $n_{i}$ be the first nonzero lattice point along the ray $\rho_{i}$
corresponding to $D_{i}$. We choose a decomposition of $N$ into a direct sum
of $\mathbb{Z}n_{1}$ with another lattice $N^{\prime}$ which is isomorphic
to $N/\rho_{1}$ (here we are abusing notation by referring to the quotient
lattice $N/\mathbb{Z}n_{1}$ as $N/\rho_{1}$).
Then we can write $n_{i}=b_{i}n_{1}+c_{i}n_{i}^{\prime}$
where $n_{i}^{\prime}\in N^{\prime}\cong N/\rho_{1}$ is indecomposable
(i.e. not of the form $k \, n_i^{\prime\prime}$
for some integer $k$ with $|k| \geq2$ and
$n_i^{\prime\prime} \in N^{\prime}$ ),
and $c_{i}$ is some nonnegative integer.
Now we choose the linear functional $u$ on $N$ such that $\left\langle
u,n_{1}\right\rangle=1$ and its restriction to $N^{\prime}$ is zero. Then
in the Chow group of $X$ we have the following relation
(see \cite[p.106]{Fulton}):
\begin{equation*}
\sum_{i=1}^{m}\left\langle u,n_{i}\right\rangle D_{i}=0.
\end{equation*}
When we restrict this relation to the toric subvariety $D_{1}$ then those
terms involving $D_{i}$ with $i>l$ will disappear because they are disjoint
from $D_{1}$, and using the formula on \cite[p. 108]{Fulton}, we have
\begin{equation*}
\bigotimes_{i=1}^{l} O_{D_{1}}
\left( \frac{\left\langle u,n_{i}\right\rangle}{c_i} D_{i}\right) =O_{D_{1}}.
\end{equation*}
Or equivalently, since $\left\langle u,n_1\right\rangle=1$ and
$\left\langle u_i,n_1\right\rangle=b_i$, we have
\begin{equation*}
O_{D_{1}}\left( -D_{1}\right)
=\bigotimes_{i=2}^{l}O_{D_{1}}\left(\frac{b_{i}}{c_i}D_{i}\right) .
\end{equation*}
Now under the identification $N^{\prime}\cong N/\rho_{1}$, the restriction
of the divisor $D_{i}$ to $D_{1}$ corresponds to the ray in $N/\rho_{1}$
spanned by $n_{i}^{\prime}$ when $2\leq i\leq l$. Therefore the piecewise
linear function $\Psi_{-D_{1}}^{D_{1}}:N_{\mathbb{R}}/\rho_{1}\rightarrow
\mathbb{R}$ is determined by $\Psi_{-D_{1}}^{D_{1}}\left( n_{i}^{\prime
}\right) =\frac{b_{i}}{c_i}$.  This implies the assertion of the proposition.
\end{proof}

\section{Acknowledgements}

The authors are grateful to Michael Davis for several very helpful comments,
proofs, references, and the permission to include them here. They also thank
Hugh Thomas for pointing out an error in an earlier version, and Xun Dong,
Paul Edelman, Kefeng Liu, William Messing, and Dennis Stanton for helpful conversations.


\begin{thebibliography}{99}
\bibitem{BabsonBilleraChan}E. Babson, L. Billera, and C. Chan,
\textit{Neighborly cubical spheres and a cubical lower bound conjecture},
\textrm{Israel J. Math.} \textbf{102} (1997), 297--315.

\bibitem{Bailey}W.N. Bailey, \textit{Generalized hypergeometric series},
Hafner, New York 1972.

\bibitem{BilleraBjorner}, \textit{Face numbers of polytopes and complexes}, in
\textrm{Handbook of Discrete and Computational Geometry}, J.E.~Goodman and
J.~O'Rourke, eds., CRC Press, Boca Raton/New York, 1997, pp. 291 -- 310.

\bibitem{BilleraLee}L.J. Billera and C.W. Lee, \textit{A proof of the
sufficiency of McMullen's conditions for $f$-vectors of simplicial convex
polytopes}, J. Combin. Theory Ser. A, \textbf{31} (1981), 237--255.

\bibitem{Bjorner-Coxetercomplexes}A. Bj\"orner, \textit{Some algebraic and
combinatorial properties of Coxeter complexes and Tits buildings},
\textrm{Adv. in Math.} \textbf{52} (1984), 173--212.

%\bibitem{BridsonHaefliger}
%M. Bridson and  A. Haefliger,
%\booktitle{Metric Spaces of Non-Positive Curvature (Grundlehren Der %Mathematischen Wissenschaften} {\bf  319}
%Soringer-Verlag, 1999. 
%

\bibitem{CharneyDavis}R. Charney and M. Davis, \textit{The Euler
characteristic of a nonpositively curved, piecewise Euclidean manifold},
\textrm{Pacific J. Math.}, \textbf{171} (1995), 117--137.

\bibitem{Coxeter}H.S.M. Coxeter, \textrm{Regular polytopes, 3rd edition},
Dover Publications, Inc., New York, 1973.

\bibitem{CoxKatz}D. Cox and S. Katz, \textrm{Mirror Symmetry and Algebraic
Geometry}, \textrm{Mathematical Surveys and Monographs} \textbf{68}, American
Mathematical Society, Providence, 1999.

\bibitem{Danilov}V.I. Danilov, \textit{The geometry of toric varieties},
\textrm{Russian Math. Surveys} \textbf{33}:2 (1978), 97--154.

\bibitem{Davis-Annals}M.W. Davis, \textit{Groups generated by reflections and
aspherical manifolds not covered by Euclidean space}, \textrm{Annals. Math.}
\textbf{117} (1983), 293--324.

\bibitem{Davis-Duke}M.W. Davis, \textit{Some aspherical manifolds},
\textrm{Duke Math. J.} \textbf{55} (1987), 105--139.

\bibitem{Davis-Handbook}M.W. Davis, \textit{Non-positive curvature and
reflection groups}, to appear in \textrm{Handbook of Geometric Topology}.

\bibitem{DavisJanuszkiewicz}M.W. Davis and T. Januszkiewicz, \textit{Convex
polytoes, Coxeter orbifolds and torus actions}, \textrm{Duke Math. J.}
\textbf{62} (1991), 417 -- 451.

\bibitem{DavisJanuszkiewiczScott}M.W. Davis, T. Januszkiewicz, and R. Scott
\textit{Nonpositive curvature of blow-ups}, \textrm{Selecta Math.} \textbf{4}
(1998), 491-547.

\bibitem{DavisOkun}M. W. Davis and B. Okun, \textit{Vanishing theorems for the
$L_{2}$-homology of right-angled Coxeter groups}, preprint, 1999.

\bibitem{DeLoera}J.A. de Loera, \textit{Triangulations of polytopes and
computational algebra}, Ph.D. dissertation, Cornell University, May 1995.

\bibitem{EdelmanReiner}P.H. Edelman and V. Reiner, \textit{H-shellings and
h-complexes}, \textrm{Adv. Math.} \textbf{106} (1994), 36--64.

\bibitem{Fulton}W. Fulton, \textrm{Introduction to toric varieties},
\textrm{Annals of Mathematics Studies} \textbf{31}, Princeton Univ. Press,
Princeton NJ, 1993.

\bibitem{Fulton-intersection-theory}W. Fulton, \textrm{Intersection theory},
Springer-Verlag, 1998.

\bibitem{Grunbaum}B. Gr\"{u}nbaum, \textrm{Convex polytopes}, \textrm{Pure and
Applied Mathematics} \textbf{16}, Interscience Publishers John Wiley \& Sons,
Inc., New York 1967

\bibitem{Guillemin}V. Guillemin, \textrm{Moment maps and combinatorial
invariants of Hamiltonian $T\sp n$-spaces}. \textrm{Progress in Mathematics},
\textbf{122}. Birkh\"{a}user Boston, Inc., Boston, MA, 1994.

\bibitem{Hartshorne}R. Hartshorne, \textrm{Algebraic geometry},
\textrm{Graduate Texts in Mathematics}, \textbf{52}. Springer-Verlag, New
York-Heidelberg, 1977.

\bibitem{Hirzebruch}F. Hirzebruch, \textrm{Topological methods in algebraic
geometry}, Springer-Verlag, Berlin, 1995.

\bibitem{Humphreys}J.E. Humphreys, \textrm{Reflection groups and Coxeter
groups}, \textrm{Cambridge Studies in Advanced Mathematics} \textbf{29},
Cambridge Univ. Press, Cambridge, 1990.

\bibitem{Kalai}G. Kalai, \textit{On low-dimensional faces that
high-dimensional polytopes must have}, \textrm{Combinatorica} \textbf{10}
(1990), 271--280.

\bibitem{Kawasaki}T. Kawasaki, \textit{The signature theorem for V-manifolds},
\textrm{Topology} \textbf{17} (1978), 75--83.

\bibitem{Kobayashi}S. Kobayashi, \textrm{Differential geometry of complex
vector bundles}, Publications of the Mathematical Society of Japan,
\textbf{15} Kan\^{o} Memorial Lectures, \textbf{5}. Princeton University
Press, Princeton NJ, 1987.

\bibitem{Lee}C. Lee, \textit{The associahedron and triangulations of the $n
$-gon}, \textrm{Europ. J. Combin.} \textbf{10} (1989), 551-560.

\bibitem{McMullen}P. McMullen, \textit{On simple polytopes}, \textrm{Invent.
Math.} \textbf{113} (1993), 419-444.

\bibitem{Milnor}J.W. Milnor and J.D. Stasheff, \textrm{Characteristic
classes}, \textrm{Annals of Mathematics Studies} \textbf{76}, Princeton
University Press, Princeton, N. J.; University of Tokyo Press, Tokyo, 1974.

\bibitem{Oda}T. Oda, \textrm{Convex bodies and algebraic geometry}
\textrm{Results in Mathematics and Related Areas} (3), \textbf{15},
Springer-Verlag, Berlin-New York, 1988.

\bibitem{Reiner}V. Reiner \textit{The distribution of descents and length in a
Coxeter group}, \textrm{Elec. J. Combin.}, \textbf{2} (1995), R25.

\bibitem{Stanley-g-theorem}R.P. Stanley, \textit{The number of faces of a
simplicial convex polytope}, \textrm{Adv. in Math.} \textbf{35} (1980), 236--238.

\bibitem{Stanley-generalized-h-vector}R.P. Stanley, \textit{Generalized
$H$-vectors, intersection cohomology of toric varieties, and related results},
\textrm{Adv. Stud. Pure Math.}, \textbf{11}, Commutative algebra and
combinatorics (Kyoto, 1985), 187--213.

\bibitem{Vinberg}E.B. Vinberg, \textit{Hyperbolic reflection groups},
\textrm{Russian Math. Surveys} \textbf{40} (1985), 31--75.

\bibitem{Ziegler}G.M. Ziegler, \textrm{Lectures on polytopes},
\textrm{Graduate Texts in Mathematics} \textbf{152}, Springer-Verlag, New
York, 1995.
\end{thebibliography}
\end{document}